\input amstex
\documentstyle{amsppt}
\input definiti
\input mathchar

\catcode`\@=11

 \def\AMSTeXfeatures{\Plainheads 
   \let\current@vert=\AMS@vert}

 \def\Plainheads{\sh@ftdiam=0.05em
   \getlabeldims
   \let\vshaftfill=\plnvsolidfill
   \let\hshaftfill=\plnhsolidfill
   \let\th@rhead=\plnrhead
   \let\th@lhead=\plnlhead
   \let\th@dnhead=\plndnhead
   \let\th@uphead=\plnuphead}
 
 \def\glet{\global\let}

 \def\LaTeXfeatures{\catcode`\@=11
   \ifx\@clnwd\undefined \nol@g
      \input ltxcode.tex \dol@g \fi
   \ltxheads \let\current@vert=\new@vert
   \providelto \catcode`\@=\active}

 \def\nol@g{\def\wlog{\edef\garbage}}
 \def\dol@g{\let\wlog=\wl@g} \let\wl@g=\wlog
 \nol@g 

 \newbox\ltobox
 \def\providelto{{\setbox\z@=
   \hbox{$\to$}\minharrlen=\wd\z@
   \global\setbox\ltobox=\hbox{$\activeat>>>$}}
   \def\lto{\mathrel{\copy\ltobox}}}

 \def\ltxheads{\sh@ftdiam=\@wholewidth
   \getlabeldims
   \let\vshaftfill= \ltxvsolidfill
   \let\hshaftfill=\ltxhsolidfill
   \let\th@rhead=\ltxrhead
   \let\th@lhead=\ltxlhead
   \let\th@dnhead=\ltxdnhead
   \let\th@uphead=\ltxuphead}
 {\catcode`\@=\active
   \gdef@#1{\csname #1\string@at\endcsname}
   \glet\activeat=@}
 \def\def@#1{\expandafter\def\csname #1@at\endcsname}

 \def@>#1>#2>{\@rrow R{#1}{#2}}
 \def@<#1<#2<{\@rrow L{#1}{#2}}
 \def@ V#1V#2V{\@rrow V{#1}{#2}}
 \def@ A#1A#2A{\@rrow A{#1}{#2}}
 \def@/#1/#2/#3/{\@rrow{#1}{#2}{#3}}
 \def@.{\ifodd\row\ifmmode\noharrow
     \else\leavevmode.\spacefactor3000 \fi
   \else\novarrow\fi}
 \def@={\ifodd\row\harrow\hequalfill{}{}%
   \else\varrow\vequalfill{}{}\fi}
 \def@:#1{\ifx=#1\harrow\deffill{}{}%
   \else\leavevmode\null:#1\fi}
 \def@|{\current@vert}
  \def\AMS@vert{\varrow\vequalfill{}{}}
  \def\new@vert#1|#2|{\ifodd\row
   \let\nextarrow\vertexvarrow
   \else\let\nextarrow\varrow\fi
   \nextarrow\vshaftfill{#1}{#2}}
 \def@-{\ifmmode\let\next\hl@ne
   \else\let\next\AMSatdash \fi \next}
  \def\hl@ne#1-#2-{\harrow\hshaftfill{#1}{#2}}
  \def\AMSatdash{\let\next\relax\leavevmode
    \def\next@{\ifx\next-%
      \def\next-{\futurelet\next\nextii@}%
     \else\def\next{\hbox{-}}\fi\next}%
    \def\nextii@{\ifx\next-\def\next-{\hbox{---}}%
      \else\def\next{\hbox{--}}\fi\next}%
    \futurelet\next\next@}
 \def@(#1){\tweenarrows{#1}}
 \def@[#1]{\setsp@n#1\relax\activeat}
 \def\fiberbox{\hbox{$\vcenter{\hr@le\hbox{\vr@le
   \kern1ex\vbox{\kern1.2ex}\vr@le}\hr@le}$}}
  \def\hr@le{\hrule height \sh@ftdiam}
  \def\vr@le{\vrule width \sh@ftdiam}
 \def@+#1+#2+#3+{\ifodd\row \harrow{#1}{#2}{#3}%
   \else \varrow{#1}{#2}{#3}\fi}


 \def\Dnarrfill{\vequalfill\Dnhe@d}
 \def\Uparrfill{\Uphe@d\vequalfill}
 
 \def\ontofill{\rtarrfill\kern-0.3em 
   \th@rhead\kern 0.3em} 

 \def\rtarrfill{\hshaftfill\th@rhead}
 \def\ltarrfill{\th@lhead\hshaftfill}
 \def\dnarrfill{\vshaftfill\th@dnhead}
 \def\uparrfill{\th@uphead\vshaftfill}
 \def\hequalfill{\plnhfill=}
 \def\deffill{:\plnhfill=}
 \def\plnvextfill#1{\setbox\z@
   \hbox{\the\textfont3 #1}%
   \dimen@=\dp\z@\advance\dimen@\ht\z@
   \copy\z@ \kern-\dimen@ 
   \cleaders\copy\z@ \vfill
   \kern-\dimen@ 
   \box\z@}
 \def\plnhfill#1{$\m@th\mkern-1.5mu\mathord#1\mkern-6mu
    \cleaders\hbox{$\mkern-2mu\mathord#1\mkern-2mu$}\hfill
    \mkern-6mu\mathord#1\mkern-1.5mu$}
 \def\vequalfill{\plnvextfill{\char'167}}
 \def\plnvsolidfill{\plnvextfill{\char'077}}
 \def\plnhsolidfill{\plnhfill-}
 \def\ltxhsolidfill{\leaders\hrule height\topofshaft depth\botofshaft
   \hfill}
 \def\ltxvsolidfill{\leaders\vrule width\sh@ftdiam\vfill}
 \def\hdashfill{\hd@sh\wd@sh
   \xleaders \hbox{\wd@sh\hd@sh\wd@sh}\hfill
   \wd@sh\hd@sh}
 \def\vdashfill{\vd@sh\wd@sh
   \xleaders \vbox{\wd@sh\vd@sh\wd@sh}\vfill
   \wd@sh\vd@sh}
 \def\dashed{\ifinmeasureCD\else
    \ifodd\row\option{\let\hshaftfill=\hdashfill}%
   \else\option{\let\vshaftfill=\vdashfill}\fi\fi}


 \newdimen\CDstrutht  \newdimen\CDstrutdp
   \CDstrutht=0.875\baselineskip
   \CDstrutdp=0.375\baselineskip
 \newdimen\CDstrutlen \CDstrutlen=\CDstrutht
   \advance\CDstrutlen by \CDstrutdp

 \def\CDstrut{\vrule
   height \ifnum\row=1 \z@\else\CDstrutht \fi
   depth \ifnum\row=\numrows \z@ \else\CDstrutdp \fi
   width\z@}

 \newdimen\CDarrsurr \CDarrsurr=0.375em
 \newdimen\CDdashlen
    \CDdashlen= 0.1875\baselineskip
 \newdimen\CDvarrlen \CDvarrlen=1.5\baselineskip
 \newdimen\minharrlen 
  \setbox\z@\hbox{$\longrightarrow$} \minharrlen=\wd\z@
 \newdimen\minCDharrlen \minCDharrlen=2.5em 
\newdimen \minc@lwd
\def\findminc@lwd{\minc@lwd=2\CDarrsurr
  \advance\minc@lwd\minCDharrlen}

 \newdimen\sh@ftdiam


 \newdimen\labelsurr \labelsurr=1.25 em

\newcount\sp@ncnt \sp@ncnt=\@ne
\newcount\sp@ncnt@ \sp@ncnt@=\@ne
\newdimen\@rrwd \newdimen\@rrdp


 \def\adjustbot#1{\option{\advance\@rrdp#1\relax}}
 
\def\pushvertex#1{\global\p@shlen#1\relax
   \global\let\maybepush=\dopush}


 \newdimen\p@shlen \p@shlen=\z@

 
 \let\maybepush=\relax
 \def\dopush{\ifinmeasureCD 
   \advance\locdimen by -\p@shlen 
   \else\advance \@rrwd by -\p@shlen \fi 
   \global\let\maybepush=\relax \global\p@shlen=\z@\relax}


 \def\span@ne{\global\sp@ncnt=\@ne\relax}
 \def\setsp@n#1#2{\global\sp@ncnt=#1\relax
   \ifx\relax#2\relax\else\global\sp@ncnt@=#2\relax\fi}

 \def\plnrhead{\llap{$\rightarrow\mkern-1.5mu$}}
 \def\plnlhead{\rlap{$\mkern-1.5mu\leftarrow$}}

 \def\clap#1{\hbox to \z@{\hss #1\hss}}

 \def\plndnhead{\hbox{\the\textfont3 \char'171}}
 \def\plnuphead{\hbox{\the\textfont3 \char'170}}
 \def\Dnhe@d{\hbox{\the\textfont3 \char'177}}
 \def\Uphe@d{\hbox{\the\textfont3 \char'176}}

 \def\ltxrhead{\raise\@xisheight
   \llap{\smash{\@linefnt\@getrarrow(1,0)}}}
 \def\ltxlhead{\raise\@xisheight
   \rlap{\@linefnt\@getlarrow(-1,0)}}
 \def\ltxuphead{\setbox\z@=\rlap{%
   \kern\@halfwidth\@linefnt\char'66}%
   \copy\z@\kern-\ht\z@}
 \def\ltxdnhead{\setbox\z@=\rlap{%
   \kern\@halfwidth\@linefnt\char'77}%
   \ht\z@=\z@\box\z@}

 \def\wd@sh{\kern0.5\CDdashlen}
 \def\hd@sh{\vrule height\topofshaft depth\botofshaft
    width\CDdashlen}
 \def\vd@sh{\hrule height\CDdashlen
   depth\z@ width\sh@ftdiam}

\def\xylist{14{3434}13{2414}12{1723}%
  23{1413}34{1153}11{0867}43{0707}%
  32{0580}21{0414}31{0291}41{0}}
\newcount\tgtcnt@
\def\find@xyargs{\dimen@=\@rrdp
  \advance\dimen@ by \CDstrutlen
  \tgtcnt@=\dimen@ \dimen@=\@rrwd 
  \divide\dimen@ by \@m 
  \divide \tgtcnt@ by \dimen@ 
  \expandafter\testxy\xylist\relax
  \unitlength=\@xarg\@rrdp
  \divide\unitlength by\@yarg\relax}
\def\testxy#1#2#3{\ifnum\tgtcnt@>#3
    \@xarg=#1\relax \@yarg=#2\relax
    \let\next=\ignorerest
  \else\let\next\testxy\fi\next}
\def\ignorerest#1\relax{\relax}

\let\scalefactor=\@ne
\def\SWarrow{\find@xyargs\vector
  (-\@xarg,-\@yarg)\scalefactor\hskip-\wd\@linechar}
\def\NWarrow{\find@xyargs\vector
  (-\@xarg,\@yarg)\scalefactor\hskip-\wd\@linechar}
\def\NEarrow{\find@xyargs\vector
  (\@xarg,\@yarg)\scalefactor}
\def\SEarrow{\find@xyargs\vector
  (\@xarg,-\@yarg)\scalefactor}
\def\rightupline{\find@xyargs\@linelen=\scalefactor
     \unitlength\@sline}
\def\rightdownline{\find@xyargs\@yarg=-\@yarg\relax
     \@linelen=\scalefactor\unitlength\@sline}

\def\Sim{\ifodd\row\setbox\z@=\hbox{$\sim$}\dimen@=\ht\z@
 \advance\dimen@ by -\@xisheight
  \vbox{\box\z@\kern-\@xisheight\kern\dimen@}%
  \else\hbox{$\wr$}\fi}

%
\def\harrow#1#2#3{\inmeasureCDtrue\findminarrwd
  {#2}{#3}{\sp@ncnt\minharrlen}\inmeasureCDfalse\span@ne
  \mathrel{\hbox{\options\hplace{#1}\ulabel{#2}\dlabel{#3}}}}

\def\noharrow{\harrow\hfill{}{}}
\def\vertexvarrow#1#2#3{\findarrdp \@rrwd=\z@ \setsp@n\@ne\@ne
  \vbox to \z@{\kern-1.2\CDstrutht
  \rlap{\options\vplace{#1}\llabel{#2}\rlabel{#3}}\vss}}

\newif\ifinmeasureCD
\def\measurelabel#1{\setbox\z@
  \hbox{$\scriptstyle#1\kern\labelsurr$}%
  \ifdim\wd\z@>\@rrwd \@rrwd=\wd\z@\fi}
\def\findminarrwd#1#2#3{\@rrwd=#3\relax
   \measurelabel{#1}\measurelabel{#2}}
\def\findCDarrwd#1#2{\@rrwd=\minCDharrlen
   \measurelabel{#1}\measurelabel{#2}%
  }

\newcount\row \row=\@ne \newcount\col \col=\@ne 
 \newcount\numrows 
\numrows=\@ne
 \newcount\numcols
\newcount\arrspan \newdimen\vrtxhalfwd  \newbox\tempbox

\def\DANABUG{\advance\col by \@ne
 \@rrwd=\minCDharrlen
  \advance\@rrwd by \vrtxhalfwd
  \advance\@rrwd by \CDarrsurr
  \ifnum\col>\numcols \numcols=\col
     \newlocdimen{col\the\col}\locdimen=\@rrwd 
  \else \ifdim\@rrwd>\c@l \c@l=\@rrwd\fi\fi}

\def\drop#1\\{
  \findvrtxhalfsum\DANABUG\advance\row by 2 \measureinit}

\def\measureinit{\col=\@ne \vrtxhalfwd=-\CDarrsurr\arrspan=\@ne\@rrwd=\z@
   \setbox\tempbox=\hbox\bgroup$}
\def\measure{
  \let\harrow\measureCDarrow
  \let\CDCR=\measureCR 
   \findminc@lwd 
  \inmeasureCDtrue
  \row=\@ne \numcols=\z@ \measureinit}

\def\endmeasure{\findvrtxhalfsum\DANABUG
  \numrows=\row 
  \inmeasureCDfalse}




\def\newlocdimen#1{\advance\dimenc@unt by \@ne
  \ifnum\dimenc@unt<\insc@unt
     \else\errmessage{No room for the CD}\fi
  \dimendef\locdimen=\dimenc@unt
  \expandafter\dimendef\csname#1\endcsname=\dimenc@unt}

 \def\r@wc@l{\csname row\the\row col\the\col\endcsname}
 \def\c@l{\csname col\the\col\endcsname}

 \def\findvrtxhalfsum{$\egroup
  \newlocdimen{row\the\row col\the\col}
  \locdimen=\vrtxhalfwd 
  \vrtxhalfwd=0.5\wd\tempbox 
  \advance\vrtxhalfwd by \CDarrsurr
  \advance\locdimen by \vrtxhalfwd 
  \advance\@rrwd by \locdimen 
  \maybepush
  \divide\@rrwd by \arrspan\relax
  \ifdim\@rrwd<\minc@lwd
    \ifnum\col>\@ne \@rrwd=\minc@lwd\fi \fi
  \loop 
    \ifnum\col>\numcols \numcols=\col
       \newlocdimen{col\the\col}
       \locdimen=\@rrwd 
    \else \ifdim\@rrwd>\c@l \c@l=\@rrwd\fi \fi
   \ifnum\arrspan>\@ne
      \advance\arrspan by -1 \advance\col by \@ne
  \repeat }

 \def\measureCDarrow#1#2#3{\findvrtxhalfsum
   \arrspan=\sp@ncnt\relax\global\sp@ncnt=1\relax
   \advance\col by \@ne
   \findCDarrwd{#2}{#3}%
   \setbox\tempbox=\hbox\bgroup$}

 \newcount\dr@tn \dr@tn=\z@
 \def\locate#1:#2{\ifinmeasureCD\else
   \count@=-#1
   \multiply\count@ by 2
   \advance\count@ by #2
   \dimen@=\count@\@rrwd
   \ifnum\dr@tn=\@ne\relax \else\dimen@=-\dimen@ \fi
   \dimen@i=\@rrdp
   \ifnum\dr@tn>\z@\advance\dimen@i by \CDstrutlen \fi
   \dimen@i=\count@\dimen@i
   \count@=#2 \multiply\count@ by 2
   \divide\dimen@ by \count@
   \divide\dimen@i by \count@
   \lift\dimen@i\nudge\dimen@\fi}

\def\betweenCDrows{\advance\row by \@ne \col=\@ne
\options}


\def\hbegin{\hbox\bgroup\kern\c@l \kern-\r@wc@l$}
\def\hend{$\glet\maybepush\relax \CDstrut\egroup}
\def\vbegin{\setbox\tempbox=\hbox\bgroup$}
\def\vend{$\egroup\ht\tempbox=\z@\dp\tempbox\CDvarrlen
  \box\tempbox}
\def\setCD{\let\harrow=\setCDarrow
  \let\CDCR=\setCR 
  \row=\@ne \col=\@ne \hbegin}
\let\endsetCD=\hend 

\def\findarrwd{\@rrwd=\z@ \count@=\col \advance\count@ by\sp@ncnt
  \loop\ifnum\count@>\col \advance\count@ by -1
      \advance\@rrwd by\csname col\the\count@\endcsname\repeat}
\def\setCDarrow#1#2#3{\kern\CDarrsurr\advance\col by \@ne
  \findarrwd \advance\@rrwd by -\r@wc@l  
  \@rrdp=\z@ 
  \maybepush
  \advance\col by -\@ne \advance\col by \sp@ncnt \span@ne
  \hbox to \@rrwd{\options
   \@rrwd=\scalefactor\@rrwd\hss
   \hplace{#1}\ulabel{#2}\dlabel{#3}\hss}%
   \kern\CDarrsurr}

\newdimen\labspacei 
\newdimen\labspaceii 

\newdimen\@xisheight
  \@xisheight=\the\fontdimen22\textfont2
\newdimen\labelskip
  \labelskip=\the\fontdimen10\textfont3 
\newdimen\topofshaft
\newdimen\botofshaft
\newdimen\botofulabel
\newdimen\topofdlabel
\def\getlabeldims{
  \topofshaft=0.5\sh@ftdiam
  \botofshaft=\topofshaft
  \advance\topofshaft by \@xisheight  
  \advance\botofshaft by -\@xisheight  
  \botofulabel=\topofshaft
  \advance\botofulabel by \labelskip
  \topofdlabel=\botofshaft
  \advance\topofdlabel by \labelskip}

\def\ulabel{\ifnum\row=\@ne\let\next\ulabeli
   \else\let\next\ulabellap\fi\next}
\def\ulabeli#1{\vbox{
  \clap{\kern-\@rrwd$\scriptstyle#1$}%
  \kern\botofulabel}\maybeoffset}
\def\ulabellap#1{\vbox to \z@{\vss
  \clap{\kern-\@rrwd$\scriptstyle#1$}%
  \kern\botofulabel}\maybeoffset}
\def\dlabel{\ifnum\row=\numrows\let\next\dlabeli
   \else\let\next\dlabellap\fi\next}
\def\dlabeli#1{\vtop{\kern\topofdlabel
  \clap{\kern-\@rrwd$\scriptstyle#1$}%
  }\maybeoffset}
\def\dlabellap#1{\vbox to \z@{\kern\topofdlabel
  \clap{\kern-\@rrwd$\scriptstyle#1$}%
  \vss}\maybeoffset}
\def\rlabel#1{\vbox to \z@{\vss
  \rlap{\kern\labelskip$\scriptstyle#1$}%
  \vss\kern-\@rrdp}\maybeoffset}
\def\llabel#1{\vbox to \z@{\vss
  \llap{$\scriptstyle#1$\kern\labelskip}%
  \vss\kern-\@rrdp}\maybeoffset}
\def\swlabel#1{\vtop{\kern0.5\@rrdp
  \llap{$\scriptstyle#1$\kern\labelskip\kern-0.5\@rrwd}
  }\maybeoffset}
\def\nwlabel#1{\vbox{
  \llap{$\scriptstyle#1$\kern\labelskip\kern-0.5\@rrwd}%
  \kern-0.5\@rrdp}\maybeoffset}
\def\selabel#1{\vtop{\kern0.5\@rrdp
  \rlap{\kern0.5\@rrwd\kern\labelskip$\scriptstyle#1$}%
  }\maybeoffset}
\def\nelabel#1{\vbox{
  \rlap{\kern0.5\@rrwd\kern\labelskip$\scriptstyle#1$}%
  \kern-0.5\@rrdp}\maybeoffset}
\def\cplace#1{\vbox to \z@{\vss
  \clap{$#1$\kern-\@rrwd}%
  \kern-\@rrdp\vss}\maybeoffset}
\def\hplace#1{\hbox to \@rrwd{#1}\maybeoffset}
\def\vplace#1{\clap{\vbox to \z@{#1\kern-\@rrdp}}\maybeoffset}

\newdimen\nudgeamount \nudgeamount=\z@
\newdimen\liftamount \liftamount=\z@
\let\maybeoffset\relax
\newbox\offsetbox \newdimen\lastheight
\def\dooffset{
  \setbox\offsetbox=\lastbox \lastheight=\ht\offsetbox 
  \setbox\offsetbox=\vbox{\kern-\liftamount\box\offsetbox}%
  \ht\offsetbox=\lastheight
  \kern\nudgeamount\box\offsetbox\kern-\nudgeamount
  \global\nudgeamount=\z@ \global\liftamount=\z@
  \glet\maybeoffset=\relax}
\def\nudge#1{\ifinmeasureCD\else
  \global\advance\nudgeamount#1\relax
  \global\let\maybeoffset\dooffset\fi}
\def\lift#1{\ifinmeasureCD\else
  \global\advance\liftamount#1\relax
  \global\let\maybeoffset\dooffset\fi}

\def\findarrdp{\@rrdp=\CDvarrlen
  \ifnum\sp@ncnt@>1
    \advance\@rrdp by \CDstrutlen
    \multiply\@rrdp by \sp@ncnt@
    \advance\@rrdp by -\CDstrutlen \fi
 }

\def\varrow#1#2#3{\ifnum\sp@ncnt>\@ne 
     \sp@ncnt@=\sp@ncnt\relax\fi
  \findarrdp \@rrwd=\z@ 
  \kern\c@l
   \hbox to \z@{\options
   \@rrdp=\scalefactor\@rrdp
    \hss\vplace{#1}\llabel{#2}\rlabel{#3}\hss}%
  \global\advance\col by \@ne \setsp@n\@ne\@ne
  }

\def\novarrow{\varrow\vfill{}{}}

\def\tweenarrows#1{\findarrwd \findarrdp \setsp@n\@ne\@ne
  \rlap{\options\cplace{#1}}}

\def\usarrow #1#2#3{\dr@tn=\@ne
  \findarrwd \findarrdp \setsp@n\@ne\@ne 
  \rlap{\options\cplace{#1}\nwlabel{#2}\selabel{#3}}%
  \dr@tn=\z@}
\def\dsarrow #1#2#3{\dr@tn=\tw@
  \findarrwd \findarrdp \setsp@n\@ne\@ne 
  \rlap{\options\cplace{#1}\swlabel{#2}\nelabel{#3}}%
  \dr@tn=\z@}
 \def\@rrow#1{\csname #1@rrow\endcsname}
 \def\R@rrow{\harrow \rtarrfill}
 \def\L@rrow{\harrow \ltarrfill}
 \def\V@rrow{\varrow \dnarrfill}
 \def\A@rrow{\varrow \uparrfill}
 \def\SE@rrow{\dsarrow \SEarrow}
 \def\NW@rrow{\dsarrow \NWarrow}
 \def\SW@rrow{\usarrow \SWarrow}
 \def\NE@rrow{\usarrow \NEarrow}
 \def\DS@rrow{\dsarrow \dnslope}
 \def\US@rrow{\usarrow \upslope}
 \def\upslope{\find@xyargs
       \@linelen=\unitlength\@sline}
 \def\dnslope{\find@xyargs\@yarg=-\@yarg\relax
       \@linelen=\unitlength\@sline}

\newtoks\optionlist 
\optionlist={}
\let\options\relax
\def\dooptions{\the\optionlist\global\optionlist={}%
  \glet\options=\relax}
\def\option#1{\ifinmeasureCD\else
  \glet\options=\dooptions
  \global\optionlist=\expandafter{\the\optionlist\relax#1}\fi}
\def\wider#1{\ifinmeasureCD\else
   \option{\advance\@rrwd by #1}\fi}
\def\deeper#1{\ifinmeasureCD\else
   \option{\advance\@rrdp by #1}\fi}


{\def\\{\global\let\sptoken= }\\ }

\def\CR{\futurelet\nexttok\testCR}
\def\testCR{\ifx\nexttok\sptoken
   \let\next\eatspaceCR\else\let\next\CDCR\fi\next}
\def\eatspaceCR#1 {\CR}
\def\measureCR{\ifx\nexttok\endmeasure\let\nextCR\relax
    \else\let\nextCR\drop\fi\nextCR}
\def\setCR{\ifodd\row
  \ifx\nexttok\endsetCD\else\hend\betweenCDrows\vbegin\fi
  \else\vend\betweenCDrows\hbegin\fi}

\countdef\dimenc@unt=11
\def\CD#1\endCD{
   \begingroup\let\\=\CR
  \m@th\offinterlineskip
   \measure#1\endmeasure\null\,\vcenter{\setCD#1\endsetCD}\,
   \endgroup
    }

\ifx\@clnwd\undefined \nol@g\else\catcode`\ =14\relax\fi
 \font\@linefnt=line10 
 \newcount\@tempcnta
 \newcount\@tempcntb
 \newdimen\@tempdima
 \newdimen\@tempdimb
 \newdimen\@wholewidth
 \newdimen\@halfwidth
   \@wholewidth\fontdimen8\@linefnt \@halfwidth .5\@wholewidth
 \newdimen\unitlength
 \newcount\@xarg
 \newcount\@yarg
 \newcount\@yyarg
 \newbox\@linechar
 \newdimen\@linelen
 \newdimen\@clnwd
 \newdimen\@clnht
 \newif\if@negarg
 
 \def\@whilenoop#1{}

 \def\@whiledim#1\do #2{\ifdim #1\relax#2\@iwhiledim{#1\relax#2}\fi}

 \def\@iwhiledim#1{\ifdim #1\let\@nextwhile=\@iwhiledim 
         \else\let\@nextwhile=\@whilenoop\fi\@nextwhile{#1}}

 \def\@sline{\ifnum\@xarg< 0 \@negargtrue \@xarg -\@xarg \@yyarg -\@yarg
   \else \@negargfalse \@yyarg \@yarg \fi
 \ifnum \@yyarg >0 \@tempcnta\@yyarg \else \@tempcnta -\@yyarg \fi
 \ifnum\@tempcnta>6 \@badlinearg\@tempcnta0 \fi
 \ifnum\@xarg>6 \@badlinearg\@xarg 1 \fi
 \setbox\@linechar\hbox{\@linefnt\@getlinechar(\@xarg,\@yyarg)}%
 \ifnum \@yarg >0 \let\@upordown\raise \@clnht\z@
    \else\let\@upordown\lower \@clnht \ht\@linechar\fi
 \@clnwd=\wd\@linechar
 \if@negarg \hskip -\wd\@linechar \def\@tempa{\hskip -2\wd\@linechar}\else
      \let\@tempa\relax \fi
 \@whiledim \@clnwd <\@linelen \do
   {\@upordown\@clnht\copy\@linechar
    \@tempa
    \advance\@clnht \ht\@linechar
    \advance\@clnwd \wd\@linechar}%
 \advance\@clnht -\ht\@linechar
 \advance\@clnwd -\wd\@linechar
 \@tempdima\@linelen\advance\@tempdima -\@clnwd
 \@tempdimb\@tempdima\advance\@tempdimb -\wd\@linechar
 \if@negarg \hskip -\@tempdimb \else \hskip \@tempdimb \fi
 \multiply\@tempdima \@m
 \@tempcnta \@tempdima \@tempdima \wd\@linechar \divide\@tempcnta \@tempdima
 \@tempdima \ht\@linechar \multiply\@tempdima \@tempcnta
 \divide\@tempdima \@m
 \advance\@clnht \@tempdima
 \ifdim \@linelen <\wd\@linechar
    \hskip \wd\@linechar
   \else\@upordown\@clnht\copy\@linechar\fi}
 
 \def\@getlinechar(#1,#2){\@tempcnta#1\relax\multiply\@tempcnta 8
 \advance\@tempcnta -9 \ifnum #2>0 \advance\@tempcnta #2\relax\else
 \advance\@tempcnta -#2\relax\advance\@tempcnta 64 \fi
 \char\@tempcnta}
 
 \def\vector(#1,#2)#3{\@xarg #1\relax \@yarg #2\relax
 \@tempcnta \ifnum\@xarg<0 -\@xarg\else\@xarg\fi
 \ifnum\@tempcnta<5\relax
 \@linelen=#3\unitlength
 \ifnum\@xarg =0 \@vvector 
   \else \ifnum\@yarg =0 \@hvector \else \@svector\fi
 \fi
 \else\@badlinearg\fi}
 
 \def\@svector{\@sline
 \@tempcnta\@yarg \ifnum\@tempcnta <0 \@tempcnta=-\@tempcnta\fi
 \ifnum\@tempcnta <5
   \hskip -\wd\@linechar
   \@upordown\@clnht \hbox{\@linefnt  \if@negarg 
   \@getlarrow(\@xarg,\@yyarg) \else \@getrarrow(\@xarg,\@yyarg) \fi}%
 \else\@badlinearg\fi}
 
 \def\@getlarrow(#1,#2){\ifnum #2 =\z@ \@tempcnta='33\else
 \@tempcnta=#1\relax\multiply\@tempcnta \sixt@@n \advance\@tempcnta
 -9 \@tempcntb=#2\relax\multiply\@tempcntb \tw@
 \ifnum \@tempcntb >0 \advance\@tempcnta \@tempcntb\relax
 \else\advance\@tempcnta -\@tempcntb\advance\@tempcnta 64
 \fi\fi\char\@tempcnta}
 
 \def\@getrarrow(#1,#2){\@tempcntb=#2\relax
 \ifnum\@tempcntb < 0 \@tempcntb=-\@tempcntb\relax\fi
 \ifcase \@tempcntb\relax \@tempcnta='55 \or 
 \ifnum #1<3 \@tempcnta=#1\relax\multiply\@tempcnta
 24 \advance\@tempcnta -6 \else \ifnum #1=3 \@tempcnta=49
 \else\@tempcnta=58 \fi\fi\or 
 \ifnum #1<3 \@tempcnta=#1\relax\multiply\@tempcnta
 24 \advance\@tempcnta -3 \else \@tempcnta=51\fi\or 
 \@tempcnta=#1\relax\multiply\@tempcnta
 \sixt@@n \advance\@tempcnta -\tw@ \else
 \@tempcnta=#1\relax\multiply\@tempcnta
 \sixt@@n \advance\@tempcnta 7 \fi\ifnum #2<0 \advance\@tempcnta 64 \fi
 \char\@tempcnta}
\catcode`\ =10

\dol@g 
\catcode`\@=\active
\LaTeXfeatures

\magnification 1200 \baselineskip 18pt


\def\pbf{\par\bigpagebreak\flushpar}
\def\pmf{\par\medpagebreak\flushpar}

\def\Ker{\hbox{\rm Ker}}
\def\Hom{\hbox{\rm Hom}}
\def\Bool{\hbox{\rm Bool}}
\def\Sub{\hbox{\rm Sub}}
\def\Val{\hbox{\rm Val}}
\def\Hal{\hbox{\rm Hal}}
\def\Ct{\hbox{\rm Ct}}
\def\Log{\hbox{\rm Log}}
\def\widetilde{\mathaccent"0365 }

\def\vp{\varphi}

\def\G{\Gamma}
\def\ol{\overline}
\def\pt{\Phi\Theta}

\def\Aut{\hbox{\rm Aut}}
\def\knowf{Know$_{\Phi \Theta}(f)$}
\def\know{Know$_{\Phi \Theta}$}
\def\em{{_{\Phi \Theta}}}

\def\vare{\varepsilon}
\def\refstyle#1{\uppercase{%
  \if#1A\relax \def\keyformat##1{[##1]\enspace\hfil}%
  \else\if#1B\relax
      \def\keyformat##1{\aftergroup\kern
           \aftergroup-\aftergroup\refindentwd}%
      \refindentwd\parindent 
  \else\if#1C\relax
      \def\keyformat##1{\hfil##1.\enspace}%
  \else\if#1D\relax 
      \def\keyformat##1{\hfil\llap{{##1}\enspace}} 
  \fi\fi\fi\fi}
} \refstyle{D}
\def\refsfont@{\tenpoint} 


\topmatter

\title
An Algebraic Approach to  Knowledge Bases Informational
Equivalence
\endtitle
\affil
Institute of Mathematics\\
Hebrew University  \\
Givat Ram, 91904, Jerusalem,  Israel\\
and\\
Department of Computer Science\\
Bar Ilan University\\
52900, Ramat Gan, Israel
\endaffil




\author
 B. Plotkin,\qquad T. Plotkin
\endauthor
\endtopmatter

\bigskip

\centerline{\smc Abstract:}

\medskip

In this paper we study the notion of knowledge from the positions
of universal algebra and algebraic logic. We consider first order
knowledge which is  based on first order logic. We define
categories  of knowledge and knowledge bases. These notions are
defined for the fixed subject of knowledge. The key notion of
informational equivalence of two knowledge bases is introduced. We
use the idea of equivalence of categories in this definition.  We
prove that for finite models there is a  clear way to determine
whether the knowledge bases are informationally equivalent.


\head Introduction
\endhead

This work stands at intersection of two areas: universal algebra
and category theory on the one hand and a field we call knowledge
science on the other. We view the latter as a science dealing with
languages of knowledge representation. It is strongly related to
universal algebra and can be considered as an area of mathematics
having motivation in computer science.

Knowledge theory and knowledge bases provide an important example
of the field where application of universal algebra and algebraic
logic is very natural, and their interacting with quite practical
problems arising in computer science is very productive. Another
examples of such interaction are given by relational database
theory, constraint satisfaction problem (\cite{BJ},\cite{JCP}),
theory of complexity, and by others.

One can speak about knowledge and a system of knowledge. As a
rule, a domain of knowledge or of a system of knowledge is fixed.
 In our approach only knowledge
that allows a formalization in some logic is considered. The logic
may be different. It is often oriented towards the corresponding
field of knowledge cf. \cite{G},\cite{L},\cite{S}.

In this paper we focus on the special situation of elementary
knowledge.


Elementary knowledge is considered to be a
 first order knowledge, i.e., the knowledge
that can be represented by the means of the First Order Logic
(FOL). The corresponding applied field (field of knowledge) is
grounded on some variety of algebras $\Theta$, which is arbitrary
but fixed. This variety $\Theta$ is considered as a knowledge
type. Its counterpart in database theory is the notion of datatype
$\Theta$.

 We also fix a set of symbols of relations $\Phi$.
The subject of knowledge is a triple $(G, \Phi, f)$, where $G$ is
an algebra in $\Theta$ and $f$ is a interpretation of the set
$\Phi$ in $G$. It is a model in the ordinary mathematical sense.
As a rule, we use  shorthand and write $f$ instead of $(G, \Phi,
f)$.  For the given $\Phi $ we denote the corresponding applied
field by $\Phi \Theta$.

FOL is also oriented on the variety $\Theta$.

We assume that every knowledge under consideration is represented
by three components: \pmf 1) {\it The description of  knowledge.}
It is a syntactical part of knowledge, written out in the language
of the given logic. The description reflects, what do we want to
know. \pmf 2) {\it The subject of  knowledge} which is an object
in the given applied field, i.e., an object for which we determine
knowledge. \pmf 3) {\it The content of  knowledge} (its
semantics).

The first two components are relatively independent, while the
third one is uniquely determined by the previous two. In the
theory under consideration, this third component has a geometrical
nature.  In some sense it is an algebraic set in an affine space.
If $T$ is a description  of knowledge and $(G, \Phi, f)$ is a
subject, then $T^f$ denotes the content of knowledge. We would
like to equip the content with its own structure, algebraic or
geometric, and to consider some aspects of such structure.

 We want to underline that there are three aspects in our approach
 to knowledge representation:  logical (for knowledge
 description), algebraic (for the subject of knowledge) and
 geometric (in the content of knowledge).
This geometry is of algebraic nature. However, the involved
algebra inherits some geometric intuition.

 Let us emphasize that logic (syntax) and geometry (semantics)
often interlace: its own geometry is possible in logic, while
logic is possible in geometry. In general, we can eliminate
geometry and reduce everything to logic. But this leads to
essential loss, namely we loose geometrical intuition which
supplements logic.

 We consider categories of elementary knowledge.  The language of
 categories in knowledge theory is a good way to organize and
 systematize primary elementary knowledge.  Morphisms in a knowledge
 category  give links between knowledge.  In particular, one
 can speak of isomorphic knowledge.  The categorical approach also
 allows us to use ideas of monada and comonada [ML].  It turns out that
 this provides some general views on enrichment and
 computation of knowledge.  Enrichment of a structure can be
 associated with a suitable monada over a category, while the
 corresponding computation is organized by comonada.
A knowledge base is related to a category of knowledge.

This paper is in a sense a continuation of \cite{PTP};  we repeat
some material to make the paper self-contained. However, there are
 certain changes in the approach to the basic notions in comparison to \cite{PTP}.
  The main one is that the definition of knowledge bases (KBs) equivalence uses
  the idea of categories
equivalence. To every KB it corresponds a database (DB)\cite{PTP}.
According to the principal result of the paper in the situation of
finite models KBs are equivalent if and only if the corresponding
databases are equivalent. This result is contained in the main
Theorem 4 of the paper.

The paper is organized as follows. We  include the material from
\cite{PTP} which is necessary for the understanding of the further
sections: the first four sections follow \cite{PTP} and provide a
background to what follows. For the details see \cite{Pl1},
\cite{Pl2},
\cite{Pl3}.
\bigskip

 \pbf
\heading 1. \ Algebra and logic
\endheading

\subheading{1.1 Multi-sorted algebra}
Keeping in mind applications, throughout the paper the term
algebra means multi-sorted, i.e., not necessarily one-sorted,
algebra.  We fix a set of sorts $\Gamma$.  In the considered
varieties $\Theta$ this set is finite, but it need not to be
finite in general.  We meet infinite $\G$ in the next section.

For every algebra $G\in \Theta$ we write
$$G = (G_i, i \in \G).
$$
The set of operations $\Omega$ is called the {\it signature of
algebras} in $\Theta$. Every symbol $\omega \in \Omega $ has a
type $\tau = \tau(\omega) = (i_1, \ldots, i_{n}; j), i , j\in \G$.
An operation of type $\tau$ is a mapping
$$
G_{i_1} \times \ldots \times G_{i_n} \to G_j.
$$
All operations of the signature $\Omega $ satisfy some set of
identities.  These identities define  the variety $\Theta$ of
$\G$-sorted $\Omega$-algebras.  Let us consider homomorphisms and
free algebras in $\Theta$ .  A homomorphism of algebras in
$\Theta$ has the form
$$
\mu=(\mu_i, i \in \G) \colon G = (G_i, i \in \G) \to G'= (G'_i, i
\in \G).
$$
Here $\mu_i \colon G_i \to G'_i$ are mappings of sets, coordinated
with operations in $\Omega$. A congruence $Ker \mu = (\Ker {\mu_
i}, i \in \G)$ is the kernel of a homomorphism $\mu$.

We consider multi-sorted sets $X = (X_i, i \in \G)$ and the
corresponding free in $\Theta$ algebras
$$
W = W(X) = (W_i, i \in \G).
$$

A set $X$ and a free algebra $W$ can be presented as  free union
of all $X_i$ and all $W_i$, respectively.

Every (multi-sorted) mapping $\mu: X \to G$ is extended up to a
homomorphism $\mu: W \to G$.  Denote the set of all such $\mu$ by
$\Hom (W, G)$.  If all $X_i$ are finite, we treat this set as an
affine space. Homomorphisms $\mu\colon W \to G$ are points of this
space.

For the given $G = (G_i, i \in \G)$ and $X=(X_i, i \in \G)$ we can
consider the set
$$
G^X = (G^{X_i}_i, i \in \G).
$$
It is the set of mappings
$$
\mu = (\mu_i, i \in \G)\colon X \to G.
$$
There is a  natural bijection $\Hom(W,G) \to G^X$. More
information about multi-sorted algebras can be found in [Pl1].

Now let us turn to the models.  Fix some set of symbols of
relations $\Phi$.  Every $\vp \in \Phi$ has its type $\tau = \tau
(\vp) = (i_1, \ldots, i_n)$.  A relation, corresponding to $\vp$,
is a subset in the Cartesian product $G_{i_1} \times\ldots\times
G_{i_n}$. Denote by $\Phi\Theta$ the class of models $(G, \Phi,
f)$, where $G \in \Theta$, and $f$ is a interpretation of the set
$\Phi$ in $G$. As for homomorphisms of models, they are
homomorphisms of the corresponding algebras which are coordinated
with relations.



\subheading{1.2 \ Logic} We consider logic in the given variety
$\Theta$. For every finite $X$, there is  a logical signature
$$
L = L_X=\{\vee,\wedge,\lnot,\exists x, \; \; x \in X\},
$$
where $X$ is $\mathop\bigcup\limits_{i\in\G} X_i$ for a finite
$\G$. We consider the set (more precisely, the $L$-algebra) of
formulas $L\Phi W$ over the free algebra $W=W(X)$.  This algebra
is an $L$-algebra of formulas of FOL over the given $\Theta$,
$\Phi$, and $X$.

First we define the atomic formulas. They are equalities of the
form $w\equiv w',$ with $w, w'\in W$ of the same sort and the
formulas $\vp (w_1, \ldots, w_n)$, where $w_i \in W,$ and all
$w_i$ are positioned according to the type $\tau=\tau(\vp)$ of the
relations $\vp$ and to the sorts.  The set of all atomic formulas
we denote by $M = M_X$.  Define $L\Phi W$ to be the absolutely
free $L_X$-algebra over $M_X$.

Let us consider another example of an $L_X$-algebra.

Given $W=W(X)$ and $G \in \Theta$, denote by $\Bool (W,G)$ the
Boolean algebra $\Sub (\Hom(W,G))$ of all subsets in $\Hom(W,G)$.
Define  the action of quantifiers in $\Bool (W,G)$.  Let $A$ be a
subset in $\Hom(W,G)$ and $x \in X_i$ be a variable of the sort
$i$.  Then $\mu\colon W\to G$ belongs to the set $\exists x A$ if
there exists $\nu \colon W \to G$ in $A$ such that $\mu (y) = \nu
(y)$ for every $y \in X$ of the sort $j, j \neq i$, and for every
$y \in X_i$, $y \neq x$.  Thus we get an $L$-algebra $\Bool (W,
G)$.

Now let us define a mapping
$$
\Val^X_f\colon  M_X \to \Bool (W,G),
$$
where $f$ is a model (the subject of knowledge), which realizes
the set $\Phi$ in the given $G$. If $w\equiv w'$ is an equality of
the sort $i$, then we set:
$$
\mu: W\to G \in \Val^X_f (w\equiv w') = \Val^X (w\equiv w')
$$
if $\mu_i(w)=\mu_i (w') $ in $G$.  Here the point $\mu$ is a
solution of the equation $w \equiv w'$.  If the formula is of the
form $\vp (w_1, \ldots, w_n)$, then
$$\mu \in \Val^X_f(\vp(w_1,\ldots, w_n))
$$
if $\vp(\mu(w_1),\ldots, \mu(w_n))$ is valid in the model $(G,
\Phi, f)$.  Here $\mu(w_j) = \mu_{i_j} (w_j)$, $i_j$ is the sort
of $w_j$.  The mapping $\Val^X_f$ is uniquely extended up to the
$L$-homomorphism
$$
\Val^X_f \colon L\Phi W\to \Bool (W,G).
$$
Thus, for every formula $u \in L\Phi W$ we defined its value
$\Val_f (u)$ in the model $(G,\Phi, f)$, which is an element in
$\Bool (W,G)$.

Every formula $u \in L\Phi W$ can be viewed as an equation in the
given model.  Then a  point $\mu\colon W\to G$ is the solution of
the ``equation" $u$ if $\mu \in \Val_f (u)$.

\subheading{1.3 Geometrical Aspect}

In the $L$-algebra of formulas $L\Phi W$, $W = W(X)$, we consider
its various subsets $T$.
 On the other hand, we consider subsets $A$  in the
affine space $\Hom (W, G)$, i.e., elements of the $L$-algebra
$\Bool (W, G)$. For each  model $(G, \Phi, f)$  and for these $T$
and $A$ we establish the following {\it Galois correspondence}
between sets of formulas in  $L$-algebra of formulas $L\Phi W$ and
sets of points in the space $\Hom (W, G)$:
$$
\eqalign{ T^f &= A = \mathop{\bigcap}\limits_{u\in T} \Val_f
(u),\cr
 A^f&=T=\{u|A\subset \Val_f (u) \}.\cr}
$$

Here $A=T^f$ is a locus of all points satisfying the formulas from
$ T$. We regard $T$ also as a system of "equations", where each
"equation" is represented by a formula $u$ from $T$. Every set $A$
of such kind is said to be an {\it algebraic set}
 (or closed set, or  algebraic variety), determined for the given model.
We define {\it knowledge} as
$$
(X, T, A, (G, \Phi, f)).
$$
Here $T$ is a {\it description of knowledge} and $(G, \Phi, f)$ is
a {\it subject of  knowledge}.  $A=T^f$ is a {\it content of
knowledge}, represented as an algebraic variety, $X$ is a {\it
place of knowledge} (the place, where the knowledge is situated).
A set $A$ can be regarded also as a relation between elements of
$G$ derived from equalities and relations of the basic set $\Phi$.
The relation $A=T^f$ belongs to the multi-sorted set
$$
G^X = \{ G^{X_i}_i, \; \; i \in \Gamma\}.
$$
A set $T$ of the form $T=A^f$ for some $A$ is called an $f$-{\it
closed set}.
 For an arbitrary $T$ we have its closure
 $
 T^{ff}=(T^f)^f$ and for every $A \subset \Hom(W, G)$ we have the
 closure $A^{ff} = (A^f)^f$.

 It is easy to understand that the following rule takes place:

\medskip
 {\it A formula $v$ belongs to the set $T^{ff}$ if and only if the
 formula
 $$(\mathop{\wedge}\limits_{u \in T} u) \to v$$ holds in the model
 $(G, \Phi, f)$.}
\medskip

 If the set $T$ is infinite then the corresponding formula is called {\it infinitary}.

 We want to study knowledge with
 different, changing ``places of knowledge" $X$.  In this case
 one should consider different $W = W(X)$, different ``spaces of
 knowledge" $\Hom(W(X), G)$, and different $L\Phi W(X)$.

 Free in $\Theta$ algebras $W(X)$ with finite $X$ are the objects
 of the category, denoted by $\Theta^0$.   Morphisms of this
 category $s\colon W(X) \to W(Y)$ are arbitrary homomorphisms of
 algebras.  The category $\Theta^0$ is a full subcategory in the
 category $\Theta$.

 We intend to build a new category related to the first order logic
 for the given $\Theta$. This category  will play for the  FOL the role
 similar to that of
 the category of free algebras $\Theta^0$ for the equational logic.
 With this end we turn from pure logic to algebraic logic.
 Such a transition will allow us to associate
 description of knowledge with its content in a more interesting way.
The sets of the type $T=A^f$ also look more natural.

 \heading
 2. Algebraic logic
 \endheading

\subheading{2.1 The main idea} Algebraic logic deals with
algebraic structures, related to
 various logical structures which correspond to different  logical calculi.
For example, Boolean algebras are associated with classical
propositional logic, Heyting algebras are associated
 with non-classical propositional logic, Tarski cylindric algebras
 and Halmos polyadic algebras are associated with FOL.

 Every logical calculus assumes that there are  formulas of the calculus,
 axioms of logic and rules of inference.  On this basis a
syntactical  equivalence of formulas compatible with their
semantical equivalence is defined.  The transition from pure logic
to algebraic logic is grounded on treating logical formulas up to
a certain equivalence. We call the corresponding classes the {\it
compressed formulas}. This transition leads to various special
algebraic structures, in particular to the structures mentioned
above.

 Every logical calculus  is usually associated with some
 infinite set of variables.  Denote such a set by $X^0$.  In our
 situation it is a multi-sorted set $X^0 = (X^0_i, i \in \Gamma)$.
 Keeping in mind theory of knowledge and its geometrical aspect we
 will use a system of all finite subsets $X = (X_i, \; i \in \G)$
 of $X^0$ instead of this infinite universum.  This gives rise to
 multi-sorted logic and multi-sorted algebraic logic.  Every formula
 has a definite type (sort) $X$.  Denote the new set of sorts by
 $\G^0$.  It is a set of all finite subsets of the initial set
 $X^0$.

\subheading{2.2 Halmos Categories} Fix some variety of algebras
$\Theta$.   This means that a
 finite set of sorts $\G$, a signature $\Omega = \Omega (\Theta)$
 related to $\G$, and a system of identities $Id(\Theta)$ are
 given.

 Define {\it Halmos categories} for the given $\Theta$.

 First, for the given Boolean algebra $B$ we define its
 existential quantifiers [HMT].
 Existential quantifiers are the mappings $\exists\colon B \to B$
 with the conditions:

 1) \ $\exists 0 = 0$,

 2) \ $a < \exists a$,

 3) \ $\exists (a \wedge \exists b) = \exists a \wedge \exists b$,
 $0, a, b \in B$.

 The universal quantifier $\forall \colon B\to B$ is defined
 dually:

 1) \ $\forall 1  = 1$,

 2) \ $a > \forall a$,

 3) \ $\forall (a\vee \forall b) = \forall a \vee \forall b.$



 Let $B$ be a Boolean algebra and $X$ a set.  We say that $B$ is a
 {\it quantifier $X$-algebra} if a quantifier $\exists
 x \colon B \to B$ is defined  for every $x \in X$
 and for every two elements $x, y \in
 X$ the equality $\exists x \exists y = \exists y \exists x$ holds.

 One may consider also {\it quantifier $X$-algebras $B$ with equalities}
over $W(X)$. In such algebras, to
 each pair of elements $w, w' \in W(X)$ of the same sort it
 corresponds an element $w \equiv w' \in B$ satisfying the
 conditions

1)  $w\equiv w$ is the unit in $B,$

2)  $(w_1\equiv w'_1 \wedge\ldots \wedge w_n \equiv w'_n) < (w_1
\ldots w_n \omega \equiv w'_1\ldots w'_n \omega) $ where $\omega$
is an operation in $\Omega$ and everything is compatible with the
type of operation.

Now we will give the general definition of the Halmos category for
the given $\Theta$, which will be followed by examples.

{\it Halmos category} $H$ for an arbitrary finite $X=(X_i, i \in
\Gamma)$ fixes some quantifier $X$-algebra $H(X)$ with equalities
over $W(X)$. $H(X)$ are the  objects  of $H$.

The morphisms in $H$ correspond to morphisms in the category
$\Theta^0$. Every morphism $s_*$ in $H$ has the form
$$
s_* = s\colon H(X)\to H(Y),
$$
where  $s \colon W(X) \to W(Y) $ is a morphism in $\Theta^0$.

We assume that

1)  The transitions $W(X) \to H(X)$ and $s \to s_*$ yield a
(covariant) functor $\Theta^0 \to H$.

2) Every $s_* \colon H(X) \to H(Y)$ is a Boolean homomorphism.

3)  The coordination with the quantifiers is as follows:

$\qquad $ 3.1)  $s_1 \exists x a = s_2 \exists x a, \quad a \in
H(X)$, if $s_1 y = s_2y$ for every $y
 \in X, \; y \neq
x$.

$\qquad $ 3.2) $s\exists x a = \exists (sx) (sa) $ if $ sx = y \in
Y$ and $y = sx$ is not in the support of $sx'$, $x' \in X, \; x'
\neq x$.

4)  The following conditions describe coordination with equalities

$\qquad $ 4.1) $s_* (w\equiv w') = (sw \equiv sw')$ for $s\colon
W(X) \to W(Y)$, $w, w' \in W(X)$ are of the same sort.

$\qquad $ 4.2) $s^x_w a \wedge (w \equiv w') < s^x_{w'} a $ for an
arbitrary $a \in H(X), x \in X, w, w'$ of the same sort with $x$
 in $W(X)$, and $s^x_w\colon W(X) \to W(X)$ is defined by the rule:
$ s^x_w (x) = w, sy = y, y \in X,\; \;  y \neq x$.

This completes the definition of the Halmos category for a given
$\Theta$.

\subheading{2.3 The example $\Hal_\Theta(G)$}

 Fix an algebra $G$ in the variety $\Theta$.  Define the Halmos
 category $\Hal_\Theta (G)$ for the given $G$.  Take a finite set
 $X$ and consider the space $\Hom(W(X), G)$.  We have defined the
 action of quantifiers $\exists x $ for all $ x\in X$ in the
 Boolean algebra $\Bool (W(X), G)$.  The equality $w\equiv w'$ in
$ \Bool
 (W(X),G)$ is defined as a diagonal, coinciding with the set of
 all $\mu: W(X) \to G$ for which $w^\mu = {w'}^{\mu}$ holds.
 It is easy to check that in this case the algebra $\Bool (W(X), G)$
 turns out to be a quantifier $X$-algebra with equalities.  We set
 $$
 \Hal_\Theta (G) (X) = \Bool (W(X), G).
 $$
 Let now $s\colon W(X) \to W(Y)$ be given in $\Theta^0$.  We have:
$$
\tilde s\colon \Hom (W(Y), G) \to \Hom (W(X), G)$$ defined by
$\tilde s (\nu) = \nu s $ for any  $\nu \colon W(Y) \to G$.

Now, if $A$ is a subset in $\Hom (W(X), G)$, then $\nu \in s_*A =
s A$ if and only if $\tilde s (\nu) = \nu s \in A$ We have a
mapping:
$$
s_* \colon \Bool (W(X), G) \to \Bool (W(Y), G)
$$ which is a Boolean homomorphism.  One can also check that $s_*$
satisfies the conditions 3--4, thereby defining the Halmos
category $\Hal_\Theta (G)$.

Note that  a conjugate mapping
$$s^*\colon \Bool (W(Y), G) \to \Bool (W(X), G),$$ where the set
$s^*B$ is the $\tilde s$-image of the set $B$ for every $B \subset
\Hom (W(Y), G)$ corresponds to each $s_*$ . Here, $s^*$ is not a
Boolean homomorphism, but it preserves sums and zero.

It may be seen that such a conjugate mapping can be defined in any
Halmos category. See, for example \cite{Pl1}.

\subheading{2.4 Multi-sorted Halmos algebras}

Fix some infinite set $X^0 = (X^0_i, i \in \G)$ and let $\G^0$ be
the set of all finite subsets $X=(X_i, i \in \G)$ in $X^0$.  In
this section multi-sorted algebra means $\G^0$-sorted.  Every such
algebra is of the form $H=(H(X), X \in \G^0)$.

 A few words about the signature of the algebras to be constructed.
First,  the signature includes $L_X$ for every $X$
 together with equalities $w\equiv w',  w, w'$ of the same  sort in
 $W(X).$ The equalities are considered as nullary operations.
  This is the signature in $H(X)$.
 Second, we consider symbols of operations of the type $s\colon W(X) \to
 W(Y)$.  To each such symbol corresponds an unary operation
$ s \colon H(X) \to H(Y)$.  Denote the
 signature consisting of all $L_X$, all equalities, and all  $s\colon W(X) \to
 W(Y)$ by $L_\Theta$.  This is the {\it signature} of FOL in
 $\Theta$ in the {\it multi-sorted variant}.

 Consider further the variety of $\Gamma^0$-sorted
 $L_\Theta$-algebras, denoted by $\Hal_\Theta$.   The identities
 of this variety exactly copy  the definition of Halmos category.
 We call algebras from $\Hal_\Theta$ {\it multi-sorted Halmos algebras}.

Every such algebra can be considered as a small Halmos category
and vice versa. Thus we come from algebra to category and back
without a special explanation.

\subheading{2.5 Algebras of formulas}

First consider a multi-sorted set of atomic formulas $M = (M(X), X
\in \G^0)$, with $M(X) = M_X$ defined as above.  All $w\equiv w'$
are viewed as symbols of nullary operations-equalities.  The set
of symbols of relations $\Phi $ is fixed.

Denote by $H_{\Phi\Theta} = (H_{\Phi \Theta} (X), X \in \G^0)$ the
absolutely free $L_\Theta$-algebra over the set $M$.  This is the
algebra of formulas of pure FOL in the given $\Theta$.

Now denote by $\tilde H_{\Phi\Theta}$ the result of factorization
of the algebra $H_{\Phi \Theta}$  by the identities of the variety
$\Hal_\Theta$.  It is the free Halmos algebra over the set of
atomic formulas $M$.

Let us introduce the following defining relations:
$$
s_* \vp (w_1,\ldots, w_n) = \vp(sw_1,\ldots, sw_n)\tag {*}
$$
for all $s\colon W(X) \to W(Y)$ and all formulas of the type
$\vp(w_1, \ldots, w_n)$ in $M(X)$.

In the sequel the principal role will play the Halmos algebra
$\Hal_\Theta(\Phi) = \Hal_{\Phi\Theta}$, defined as a quotient
algebra of the free algebra $\tilde H_{\Phi\Theta}$ by the
relations of the (*)
 type.  Elements of this algebra are defined to be {\it compressed formulas}.

 Consider now values of formulas.  First of all take a mapping
$$
\Val_f = (\Val_f^X, X \in \G^0)\colon \ M \to \Hal_\Theta (G).
$$
For the model $(G, \Phi, f)$ the mapping $\Val^X_f\colon M_X \to
\Bool (W(X), G)=\Hal_\Theta(G)(X)$ has been defined.

This mapping is uniquely extended up to the homomorphisms
$$
\eqalign{ &\Val_f\colon H_{\Phi\Theta}\to\Hal_\Theta (G),\cr
&\Val_f\colon \tilde H_{\Phi\Theta}\to\Hal_\Theta(G).\cr}
$$
Note that the relations $(*)$ hold in every algebra
$\Hal_\Theta(G)$ and  this gives a canonical homomorphism of
Halmos algebras $$Val_f :\Hal_\Theta (\Phi) \to \Hal_\Theta (G).$$
It determines the value of the formulas $\Val_f(u)$ (pure and
compressed) in the given model $(G, \Phi, f)$.

We call two pure formulas $u$ and $v$ of the given type $X$ {\it
semantically equivalent}, if $\Val_f (u) = \Val_f (v)$ for every
model $(G, \Phi, f)$.

 The following main theorem takes place \cite{Pl2}:

 \proclaim{Theorem 1}  Two formulas $u$ and $v$ are semantically
 equivalent if and only if the corresponding compressed formulas $\ol u $ and
 $\ol v$ coincide in the algebra $\Hal_\Theta (\Phi)$.
\endproclaim

This theorem explains the role of algebra $\Hal_\Theta(\Phi)$ as a
main structure of the multi-sorted algebraic logic for FOL in the
given $\Theta$.  The same algebra plays an essential part in the
algebraic geometry in the FOL  in $\Theta$.  In particular, the
role of the algebras $\Hal_\Theta (G)$ is underlined by the
following theorem \cite{Pl2}:

\proclaim{Theorem 2}  The algebras $\Hal_\Theta (G)$ over
different $G \in \Theta$ generate the variety of Halmos algebras
$\Hal_\Theta$.
\endproclaim

Define the notion of the logical kernel of a homomorphism.

Let the homomorphism $\mu\colon W(X) \to G$ be given.  One can
view its kernel $\Ker \mu$ as a system of all formulas $w \equiv
w'$ with $w, w'$ of the same sort in $W(X)$, for which $\mu \in
\Val (w \equiv w')$.

{\it Logical kernel} $\Log\Ker \mu$ naturally generates the
standard $\Ker \mu$. We set: the formula $u \in
\Hal_{\Phi\Theta}(X)$ belongs to $\Log\Ker (\mu)$ if the point
$\mu$ lies in $\Val_f (u)$, i.e.,  if $\mu$ is a solution of the
``equation" $u$ in the given model $(G, \Phi, f)$.  It is easy to
understand, that for every point $\mu$ its logical kernel is an
ultrafilter of the Boolean algebra $\Hal_{\Phi\Theta} (X)$.  It is
also clear, that the kernel $\Ker \mu$ is the set of all
equalities in the logical kernel.

\heading 3. Categories of algebraic sets
\endheading

\subheading{3.1 Preliminary remarks}

We defined in the subsection 1.3 the algebraic sets  determined by
FOL formulas. Now  we work with the compressed formulas, i.e., the
formulas of the algebra $\Hal_\Theta (\Phi) = \Hal_{\Phi\Theta}$.
Correspondingly, we have to extend the definition of Galois
correspondence from 1.3 to the case of compressed formulas, i.e,
to the elements of $\Hal_\Theta (\Phi)$.

For the given place $X$ consider sets of formulas $T$ in
$\Hal_{\Phi\Theta} (X)$ and the sets of points $A$ in the space
$\Hom (W(X), G)$. Having the  model $(G, \Phi, f)$, we establish a
Galois correspondence between  sets of elements (compressed
formulas) in the Halmos algebra $\Hal_\Theta (\Phi)$ and sets of
points in the space $\Hom (W, G)$: :
$$
\eqalign{ &T^f = A = \bigcap_{u\in T} \Val_f(u) = \{ \mu \vert
T\subset \Log\Ker (\mu)\}\cr &A^f = T = \{ u \vert A\subset \Val_f
(u) \} = \bigcap_{\mu \in A} \Log\Ker (\mu).\cr}
$$
As in 1.3, we call a set $A$ represented as $A=T^f$ an {\it
algebraic set} or {\it algebraic variety} for the given model $(G,
\Phi, f)$.

The set $T$, represented as $A^f=T$, is always a filter of the
Boolean algebra $\Hal_{\Phi\Theta} (X)$, since by  definition it
is an intersection of ultrafilters.  We call it an {\it $f$-closed
filter}. One can consider a Boolean algebra $\Hal_{\Phi\Theta}
(X)/T$ for this $T$.  If $T^f=A$ and $A^f=T$, then the algebra
$\Hal_{\Phi\Theta} (X)/T$ is considered as an invariant of the
algebraic set $A$.
 This invariant is a {\it coordinate algebra} of the set $A$.  It can be
 viewed as an algebra of regular functions determined on the
 variety $A$ (see [Pl2]).

 Suppose an algebraic set $A$ is given. A filter $T=A^f$ can be treated as the
theory of a the set $A$ for  the fixed model $(G,\Phi,f)$.

 Every algebraic set, defined in Subsection 1.3, is also an algebraic
 set according to this new definition.  The opposite is not true,
 because in the new variant additional operations of the type
 $s\colon W(X) \to W(Y)$ are involved in the formulas.

 We will return later to the structure  of algebraic
 sets.

 Consider now the relation between the Galois correspondence  and
 morphisms of Halmos categories.

 For every $ s\colon W(X) \to W(Y)$ and every $A$ of the type $X$
 we considered a set $B=s_* A$ of the type $Y$.  If $B$ is of the
 type $Y$, then $A = s^* B $ is of the type $X$.  Define the
 operations $s_*$ and $s^*$ on the sets of formulas.

 If $T$ is a set of formulas in $\Hal_{\Phi\Theta}(Y)$, then
 $s_*T$ is a set of formulas in $\Hal_{\Phi\Theta} (X)$ defined by
 the rule:
 $$
 u \in s_*T \Leftrightarrow s  u \in T.
 $$
 If $T$ is a set of formulas in $\Hal_{\Phi\Theta}(X)$, then $s^*
 T$ is contained in $\Hal_{\Phi\Theta} (Y)$ and it is defined by
 $$
 u\in s^*T \; \; \hbox{\rm if } \; \; u = sv, \; \;  \; v \in T.
 $$
 The following theorem \cite{Pl2} holds:
 \proclaim{Theorem 3}

 1. If $T$ lies in $\Hal_{\Phi\Theta} (X)$, then
 $$
 (s^* T)^f = s_* T^f = sT^f.
 $$

 2.  If $B \subset \Hom (W(Y), G)$, then
 $$
 (s^*B)^f = s_* B^f.
 $$

 3.  If $A \subset \Hom(W(X), G)$, then $s^*A^f\subset(s_*A)^f$.
\endproclaim

 It follows from these rules that

{\it
 1.  If $A = T^f$ is an algebraic set, then $sA$ is also an
 algebraic set.

 2.  If $T=B^f$ is $f$-closed, then $sT=s_*T$ is $f$-closed.
}

\subheading{3.2. Categories $K_{\Phi\Theta} (f)$ and
$C_{\Phi\Theta}(f)$}

 Fix a model $(G,\Phi, f)$ and define a category of algebraic
 sets $K_{\Phi\Theta} (f)$ for this model.  Objects of this
 category have the form $(X, A)$, where $A = T^f$ for some $T$.
 $X$ is the place for both $A$ and $T$.

 Let us now define  morphisms
 $
 (X,A)\to (Y,B).
 $
 For $s\colon W(Y) \to W(X)$ we say that $s$ is
 {\it admissible} for
$A$ and $B$ if $\tilde s(\nu) = \nu s \in B$ for any $\nu \in A$.
It is clear that $s$ is admissible for $A$ and $B$ if $A\subset
sB$. A mapping $[s]:A\to B$ corresponds to each $s$ admissible for
$A$ and $B$. Note that for the equal $[s_1]$ and $[s_2]$ the
corresponding $\tilde s_1$ and $\tilde s_2$ can be different.

We consider {\it weak} and {\it exact} categories
$K_{\Phi\Theta}(f)$. In the first one the morphisms are of the
form $\tilde s:(X,A)\to (Y,B)$, while in the second one they are
of the form $[s]:(X,A)\to (Y,B)$. Here, $s$ assumed to be
admissible for $A$ and $B$.

If $s_1$ is admissible for $A$ and $B$ and $s_2$ for $B$ and $C$,
then $A \subset s_1 B$,
 $B \subset s_2 C$, $ s_1 B \subset s_1s_2 C$, and $s_1s_2$ is
admissible for $A$ and $C$.

Define now a category $C_{\Phi\Theta} (f)$.  Its objects are
Boolean algebras of the type
\newline
 $\Hal_{\Phi\Theta} (X) / T$, where $T
= A^f$ for some $A$.

Consider morphisms
$$
\Hal_{\Phi \Theta}(Y) / T_2 \overset \ol s\to{\longrightarrow}
\Hal_{\Phi\Theta} (X) /T_1.
$$
We proceed here from $s\colon W(Y) \to W(X)$ and pass to the new
$s : \Hal_{\Phi\Theta} (Y) \to \Hal_{\Phi\Theta} (X)$. Assume that
$s u \in T_1$ for every $u \in T_2$. The homomorphism $s$ is
admissible for $T_2$ and $T_1$ in this sense.  Define
homomorphisms  $\ol s$ for such $s$.  This defines morphisms in
$C_{\Phi \Theta}(f)$.

The next two straightforward propositions determine the
correspondence between the categories $K_{\Phi\Theta}(f)$ and
$C_{\Phi\Theta} (f)$.

\proclaim{Proposition 1} A homomorphism $s\colon W(Y)\to W(X)$ is
admissible for the sets $(X, A)$ and $(Y, B)$ if and only if it is
admissible for $T_2 = B^f$ and $T_1 = A^f$.
\endproclaim
\proclaim{Proposition 2}   If $s_1, s_2 \colon W(Y) \to W(X)$ are
admissible for $A$ and $B$, then $[s_1]=[s_2]$ implies $\ol{s_1} =
\ol{s_2}$.
\endproclaim

It follows from these two propositions that the transition
$$
(X,A) \to\Hal_{\Phi\Theta} (X)/A^f
$$
determines a contravariant functor
$$
K_{\Phi\Theta} (f) \to C_{\Phi\Theta} (f)
$$
for weak and exact categories $K_{\Phi\Theta} (f)$. Duality for
these categories takes place under some additional conditions.

\subheading{3.3 Categories $K_{\Phi\Theta} $ and $C_{\Phi\Theta}$}

In the categories $K_{\Phi\Theta}$ and $C_{\Phi\Theta}$ the model
$(G, \Phi, f)$ is not fixed.  Objects of $K_{\Phi\Theta}$ have the
form $(X,A; G,f)$.  Here $f$ is a interpretation of the set $\Phi$
in the algebra $G$, fixed for the category $K_{\Phi\Theta}$, and
$A=T^f$ for some $T\subset \Hal_{\Phi\Theta}(X)$.

Define morphisms
$$
(X,A; G_1, f_1) \to (Y,B; G_2, f_2).
$$
They act on all components of the objects.  Proceed from the
commutative diagram
$$
\CD
W(Y) @>s>> W(X) \\
@V\nu' VV @VV\nu V\\
G_2 @<\delta<< G_1\\
\endCD
$$
Consider a pair $(s, \delta)$ and write $(s, \delta) (\nu) = \nu'
= \delta\nu s$.

Let now $A = T_1^{f_1} $ be of the type $X$ and $B=T_2^{f_2}$ of
the type $Y$.  We say that the pair $(s, \delta)$ is {\it
admissible} for $A$ and $B$ if $(s, \delta) (\nu) \in B$ for every
$\nu \in A$.

We need some further auxiliary remarks. For every $\delta: G_1 \to
G_2$ and every $X$ we have a mapping
$$
\tilde\delta: \Hom(W(X), G_1) \to \Hom (W(X), G_2)
$$
defined by the rule
$$
\tilde \delta(\nu) = \delta\nu, \nu \in \Hom (W(X), G_1).
$$

Define $\delta_* A\subset \Hom (W(X), G_1)$ for every $A\subset
\Hom (W(X), G_2)$, by setting
$$
\nu \in \delta_* A \; \; \hbox{\rm if } \; \; \delta\nu=\tilde
\delta(\nu) \in A.
$$
We write also  $\delta_* A=\delta A, $ and consider $\delta^*$
determined by:
if $A \subset \Hom (W(X), G_1)$, then $\delta^* A \subset \Hom
(W(X), G_2)$ and $\nu \in \delta^* A$ if $\nu = \delta\nu_1, \nu_1
\in A$.

Now we can  say that the pair $(s, \delta)$ is admissible for $A$
and $B$ if $\delta^*A \subset s B$, or,  the same, $A\subset
\delta sB = s\delta B$.

We have morphisms
$$
(s, \delta) : (X, A; G_1, f_1)\to(Y, B; G_2, f_2)$$ and
$$
([s], \delta)\colon (X, A; G_1, f_1) \to (Y, B; G_2, f_2)
$$
for the admissible $(s, \delta)$.   Here $[s]: A \to B$ is a
mapping, induced by the pair $(s, \delta)$.  We get weak and exact
categories $K_{\pt}$.  It can be proven that the pair $(s,
\delta)$ is admissible for $A$ and $B$ if and only if the
homomorphism $s\colon \Hal_{\pt} (Y)\to Hal_{\pt} (X) $ is
admissible in respect to $T_2 = B^{f_2}$ and $T_1 = (\delta^*
A)^{f_2}$.  This leads to a natural definition of the category
$C_{\pt}$ with contravariant functor $K_{\pt} \to C_{\pt}$.

Let us define the categories $K_{\pt}(G)$ and $C_{\pt} (G)$. Here
$G$ is a fixed algebra in $\Theta$, while the interpretations $f$
of the set $\Phi$ in $G$ change.

The objects in $K_{\pt}(G)$ have the form
$$
(X, A; f).
$$
The morphisms
$$
(X, A, f_1) \to (Y, B, f_2)
$$
are defined according to the general definition of the morphisms
in $K_{\pt}$ with identical $\delta = \vare \colon G\to G$.

Objects in $C_{\pt} (G)$ have the form
$$
(\Hal_{\pt} (X) /T, f), \; \; \hbox{\rm where} \; \; T=A^f
$$
for some $A$ of the type $X$.

The transition
$$
(X, A; f) \to (\Hal_{\pt} (X) /A^f, f)
$$
determines the functor $K_{\pt} (G) \to C_{\pt} (G)$. Here
$K_{\pt} (G$) is a subcategory in $K_{\pt}$ and every $K_{\pt}(f)
$ is a subcategory in $K_{\pt} (G)$.  The same holds for $C$. See
also \cite{NP}.

\heading 4. Categories of elementary knowledge
\endheading

\subheading{4.1 The category Know$_{\Phi \Theta} (f)$}

In Subsection 1.3 we defined  knowledge as
$$
(X, T, A, (G, \Phi, f)),
$$
where each component has the corresponding meaning.
 Fix a model (subject of knowledge) $(G, \Phi, f)$. Let us
define a category of knowledge for this model and denote it by
Know$_{\Phi \Theta}(f)$. This is the knowledge category for the
given subject of knowledge. Since the model is fixed, the objects
of the category \knowf \ have to have the form $(X, T, A)$.  We do
not fix the subject of knowledge in the notation of the object,
since it is fixed in the notation of the category.

  The set $X$ is
multi-sorted.  It marks the ``place" where the knowledge is
situated. The set $X$ points also  the ``place of the knowledge'',
i.e., the space of the knowledge $Hom (W(X), G)$, while the
subject of the knowledge $(G, \Phi, f)$ is given. The set $T$ is
the description of the knowledge in the algebra $Hal_{\empty
\Theta}(X)$, and $A=T^f$ is the content of knowledge, depending on
$T$ and $f$. The set $T^{ff} = A^f$ is the full description of the
knowledge $(X, T, A)$ which is a Boolean filter in Hal$_{\Phi
\Theta}(X)$.

Now about morphisms $(X, T_1, A) \to (Y, T_2, B)$. Take $s\colon
W(Y) \to W(X)$.  We have also $s\colon \Hal_{\Phi \Theta}(Y)\to
\Hal_{\Phi \Theta}(X)$ (see 2.2). This is a homomorphism of
Boolean algebras.  The homomorphism $s$ gives rise  to
$$
\tilde s \colon Hom (W(X), G) \to Hom (W(Y), G).$$ As above, the
first $s$ is admissible for $A$ and $B$ if $\tilde s (\nu) =\nu s
\in B$ for every point $\nu\colon W(X)\to G$ in $A$.

As we know, $s$ is admissible for $A$ and $B$ if and only if $su
\in A^f$ for every $u \in B^f$.  This holds for  $s_\ast$, for
which we have also a homomorphism $\ol s:\Hal\em (Y)/B^f\to
\Hal\em (X)/A^f$.  It is easy to prove that $s$ is admissible for
$A$ and $B$ if and only if $s u \in A^f$ holds for every $u \in
T_2$. We consider admissible $s$ as a morphism
$$
s\colon (X, T_1, A) \to (Y, T_2, B),
$$
in the weak category \knowf.

  We have $\tilde s (\nu) = \nu s \in
B$ if $\nu \in A$, and $s$ induces a mapping $[s]\colon A \to B$.
Simultaneously, there is a mapping $s\colon T_2 \to A^f$ and a
homomorphism
$$
\ol s\colon \Hal\em (Y)/B^f \to \Hal\em(X)/A^f.
$$
We have already mentioned (Proposition 2) that $\ol s_1 = \ol s_2$
follows from $[s_1]=[s_2]$. Thus, we can take the morphisms of the
form
$$
[s]\colon (X, T_1, A) \to (Y, T_2, B),
$$
for the morphisms of the exact category \knowf. The canonical
functors \knowf$\to K\em(f)$ for weak and exact categories are
given by the transition $(X, T, A) \to (X, A)$. In this transition
we ``forget" to fix the description of knowledge $T$.

\subheading{4.2 The category \know}

Let us define the category of elementary knowledge for the whole
applied field $\Phi\Theta$; the subject of the knowledge $(G,
\Phi, f)$ is not fixed.  As earlier, we proceed from the category
$\Phi \Theta$ whose morphisms are homomorphisms in $\Theta$. They
ignore the relations from $\Phi$.

An object of the knowledge category \know \   has the form
$$
(X, T, A; (G, \Phi, f)),
$$
and we write $(X, T, A; G,  f)$, because $\Phi$ is fixed for the
category. Here $X$ marks the place of knowledge. The components $A
= T^f$, $G$ and $f$ may change.

Consider morphisms:
$$
(X, T_1, A; G_1, f_1) \to (Y, T_2, B; G_2, f_2).
$$
We apply the same approach as in Section 3.3 with some
modifications.

Start from $s: W(Y) \to W(X)$ and $\delta: G_1 \to G_2$. These $s$
and $\delta$ should correlate. Let us explain the correlation
condition.  Take a set $A_1 = \{ \delta \nu, \nu \in A\}=\delta^*
A$ and
 take further $T^\delta_1 = A_1^{f_2}$.
 Correlation of $s$ and $\delta$ means that $su\in T^\delta_1$ holds for
 any $u \in T_2$.  The same holds for every $u \in B^{f_2}$.
 The last also says that there is a homomorphism
 $$
 \ol s\colon \Hal\em (Y) /B^{f_2} \to \Hal\em (X)/A_1^{f_2}.
 $$
 The first of the two mappings $(s, \delta) \colon A \to B$ and
 $s\colon T_2 \to T^\delta_1$ transforms the content of knowledge,
 while the second one acts on the description.  Here $T_2$ and
 $T^\delta_1$ describe knowledge associated with the same subject
 $(G_2, \Phi, f_2)$.

 With the fixed $\delta$ there is also an exact mapping $([s], \delta): A \to
 B$.  This brings us to  weak and exact categories \know.  The morphisms of
 the first one are $(s, \delta)$ and in the second one they are of the form $([s],
 \delta)$ for $(X, T_1, A; G_2, f_1) \to (Y, T_2, B;
 G_2, f_2)$.
 The canonical functors \know$ \to K\em$ are defined by the
 transition
 $$
 (X, T, A; G, f) \to (X, A; G, f).
 $$
 As above, we remove the description of knowledge from the
 notations.

\subheading{4.3 Categories $K\em(G)$ and \know$(G)$}

An algebra $G \in \Theta$ is fixed in the categories $K\em (G)$
and \know$(G)$.  A set of symbols of relations $\Phi $ is fixed as
usual, but interpretations $f$ of $\Phi$ in $G$ may change. Thus,
$K\em(G)$ is a subcategory in $K\em$ and \know$(G)$ is a
subcategory in \know.   Here the corresponding $\delta: G \to G$
are identical homomorphisms.  Objects of the category $K\em(G)$
have the form $(X, A, f)$, and those of the category \know$(G)$
are written as $(X, T, A, f)$.  There is a canonical functor
\know$(G) \to K\em(G)$. As for morphisms
$$
\eqalign{ &(X, A, f_1) \to (Y, B, f_2) \; \; \hbox{\rm and} \cr
&(X, T_1, A, f_1) \to (Y, T_2, B, f_2),\cr}
$$
we  note that $A=A_1, A_1^{f_2} = T^\delta_1$ and
$A^{f_2}=T_1^{f_1f_2}$. Hence, the corresponding admissible $s:
W(Y) \to W(X)$ transfers each $u \in T_2$ into $s u \in T_1^{f_1
f_2}$ and it induces a homomorphism
$$
\ol s: \Hal\em (Y) /B^{f_2} \to \Hal\em (X) /A^{f_2}.
$$
Every $s$ gives a mapping $[s]: A \to B$. This defines a morphism
$(X,A,f_1)\to (Y,B,f_2)$.


\heading{ 5. Knowledge bases}
\endheading

\subheading{5.1.Category of knowledge description}

Denote the category of knowledge description by $L_{\Phi\Theta}$
or $L_\Theta(\Phi)$.

Its objects are of the form $(X,T)$, where $X$ is a finite set and
$T$ is a set of formulas of $\Hal_{\Phi\Theta}(X)$. Define
morphisms $(X,T_1)\to(X,T_2)$. According to the definition of the
category
 $\Hal_{\Theta}(\Phi)$ proceed from the functor $\Theta^0\to
\Hal_{\Theta}(\Phi)$ which assigns a mapping $s_\ast:
\Hal_{\Phi\Theta}(X) \to \Hal_{\Phi\Theta}(Y)$ to every
homomorphism $s: W(X)\to W(Y)$. We say that $s$ is {\it admissible
in respect to $T_1$ and $T_2$} if $s_\ast(u)\in T_2$ for every
$u\in T_1$. For such admissible $s$ we have a mapping $s_\ast:
T_1\to T_2$ which determines
$$
s_\ast: (X, T_1)\to (X,T_2).
$$

\subheading{5.2 Functor of transition from knowledge description
to knowledge content}

Proceed from the model $(G,\Phi, f)$ and consider a functor
$$
\Ct_f:L_{\Phi\Theta}\to K_{\Phi\Theta}(f).
$$
Here, $K_{\Phi\Theta}(f)$ is the corresponding category  of
algebraic (elementary) sets over the given model and $\Ct$ stands
for "contents". The functor $\Ct_f$ is a contravariant one. To
every object $(X,T)$ of the category $L_{\Phi\Theta}$ it assigns
the corresponding content $(X, T^f)=(X,A)$ which is an object of
the category $K_{\Phi\Theta}(f)$.

Now one has to define the functor $\Ct_f$ on morphisms. Let  a
morphism
$$
s_\ast: (Y, T_2)\to (X,T_1)
$$
be given for $s: W(Y)\to W(X)$. Show that $s$ induces a morphism
$$
\widetilde {s_\ast}: (X, A)\to (Y,B),
$$
where $A=T_1^f$, and $B=T_2^f$.

We proceed from $\widetilde s: \Hom(W(X),G)\to \Hom (W(Y),G)$.

Let us define a transition $s\to \widetilde s.$ 

Check first that if $s$ is admissible for $T_2$ and $T_1$ then
this $ s$ is admissible for $A=T_1^f$ and $B=T_2^f$. The last
means that $\widetilde s(\nu)\in B$ if $\nu\in A$. The inclusion
$\nu\in A$ says that $\nu\in \Val_f(v)$ for every $v\in T_1$. We
need to verify that $\nu s\in B,$ that is $\nu s \in \Val_f(u)$
for every $u\in T_2$.

Take an arbitrary $u\in T_2$. We have: $v= {s_\ast}(u)\in T_1$;
$\nu\in \Val_f(v)=\Val_f(s_\ast u)= s \Val_f(u)$. This gives $\nu
s  \in \Val_f(u)$. We used that $s$ and $\Val_f$ commute, since $
\Val_f$ is a homomorphism of algebras.

 The mapping $[s]: A\to B$
corresponds to the homomorphism $s:W(Y)\to W(X)$. This mapping is
considered simultaneously as a morphism in the category
$K_{\Phi\Theta}(f)$ (see 3.2)
$$
[s]: (X,A)\to (Y,B).
$$
We define: $\Ct_f(s_\ast)=\widetilde s_\ast=[s]$.

 Check now compatibility of the definition of
$\Ct_f$ with the multiplication of morphisms. Given $s_1: W(X)\to
W(Y)$ and $s_2: W(Y)\to W(Z)$ we have $s_2s_1: W(X)\to W(Z)$.
Using the fact that the transition $\Theta^0\to \Hal_\Theta(\Phi)$
is a functor, we get $(s_2s_1)_\ast= s_{2\ast}s_{1\ast}$. Here, we
have
$$
s_{1\ast}:\Hal_{\Phi\Theta}(X)\to \Hal_{\Phi\Theta}(Y),
$$
$$
s_{2\ast}:\Hal_{\Phi\Theta}(Y)\to \Hal_{\Phi\Theta}(Z),
$$
and
$$
(s_2s_1)_\ast:\Hal_{\Phi\Theta}(X)\to \Hal_{\Phi\Theta}(Z).
$$
Let $(X,T_1)$, $(Y,T_2)$ and $(Z,T_3)$ be objects in
$L_\Theta(\Phi)$ , and $s_1,s_2$ admissible in respect to $T_1$,
$T_2$ and, correspondingly, for $T_2$, $T_3$. In this case there
are morphisms
$$
s_{1\ast}:(X, T_1)\to (Y,T_2),
$$
$$
s_{2\ast}:(Y, T_2)\to (Z,T_3),
$$
and
 $$ s_{2\ast}s_{1\ast}= (s_2s_1)_\ast: (X,T_1)\to (Z,T_3).$$

Take $T_1^f=A$, $T_2^f=B$, $T_3^f=C$. We have
$$
\widetilde {s_{1\ast}}: (Y,B)\to (X,A),
$$
$$
\widetilde {s_{2\ast}}: (Z,C)\to (Y,B),
$$
and
$$
\widetilde{{s_2s_1}_\ast}=\widetilde{s_{1\ast}}\widetilde{s_{2\ast}}:
(Z,C)\to (X,A).
$$
This gives compatibility of the functor $\Ct_f$ with the
multiplication of morphisms. Compatibility with the unity morphism
is evident. This finishes the definition of the contravariant
functor $\Ct_f: L_{\Phi\Theta}\to K_{\Phi\Theta}(f)$.

\subheading{ 5.3 Homomorphisms of Halmos algebras
$\Hal_\Theta(\Phi)$ and functors of the categories
$L_\Theta(\Phi)$}

Given a homomorphism $\beta:\Hal_\Theta(\Phi_1)\to
\Hal_\Theta(\Phi_2)$,  define the corresponding functor
$\widetilde\beta:L_\Theta(\Phi_1)\to L_\Theta(\Phi_2)$. For every
set of formulas $T\subset \Hal_{\Phi_1\Theta}(X),$ denote by
$T^\beta$ the set $T^\beta=\{u^\beta, u\in T\}$. If $(X,T)$ is an
object in $L_\Theta(\Phi_1)$, then, setting
$$
\widetilde\beta(X,T)=(X,T^\beta),
$$
we get an object in $L_\Theta(\Phi_2)$.

In order to define the functor $\widetilde\beta$ on morphisms let
us make a remark. Proceed from the functors
$\Theta^0\to\Hal_\Theta(\Phi_1)$ and
$\Theta^0\to\Hal_\Theta(\Phi_2)$. The morphisms
$$
s_{\ast}^1:\Hal_{\Phi_1\Theta}(X)\to \Hal_{\Phi_1\Theta}(Y),
$$
$$
s_{\ast}^2:\Hal_{\Phi_2\Theta}(X)\to \Hal_{\Phi_2\Theta}(Y)
$$
correspond to every $s: W(X)\to W(Y)$. We have also
$$
\beta=(\beta_X, X\in \Gamma^0):\Hal_{\Theta}(\Phi_1)\to
\Hal_{\Theta}(\Phi_2).
$$
The fact that the homomorphism $\beta$ is compatible with the
operation $s$ is represented by the commutative diagram

$$
\CD
 \Hal_{\Phi_1\Theta}(X)@>s_\ast^1>> \Hal_{\Phi_1\Theta}(Y)\\
@V\beta_X VV @VV\beta_Y V\\
\Hal_{\Phi_2\Theta}(X)@>s_\ast^2>> \Hal_{\Phi_2\Theta}(Y)\\
\endCD
$$
So, for a homomorphism $s: W(X)\to W(Y)$ we have the equality
$\beta_Ys_\ast^1(u)=s_\ast^2\beta_X(u)$ for every $u\in
\Hal_{\Phi_1\Theta}(X).$

Now we are able to define an action of the functor
$\widetilde\beta$ on morphisms. Let a morphism $s_\ast^1:
(X,T_1)\to(Y,T_2)$ in the category $L_{\Phi_1\Theta}$ be given and
$s_\ast^1(u)\in T_2$ if $u\in T_1$. Then, we have $s_\ast^2(v)\in
T_2^\beta$ if $v\in T_1^\beta$.

Indeed, let $v=\beta_X(u)$, $u\in T_1$, $v\in T_1^{\beta_X}$. We
have:
$$
s_\ast^2\beta_X(u)=s_\ast^2(v)=\beta_Ys_\ast^1(u)\in
T_2^{\beta_Y},
$$
since $s_\ast^1(u)\in T_2$. Hence, $s_\ast^2(v)\in T_2^{\beta_Y}$
for every $v=\beta_X(u)\in T_1^{\beta_X}$.

 We set
$s_\ast^2=\widetilde\beta(s_\ast^1): T_1^{\beta_X}\to
T_2^{\beta_Y}$.
 A morphism
 $$
 s_\ast^2=\widetilde\beta(s_\ast^1):
(X,T_1^{\beta_X})\to (Y,T_2^{\beta_Y})
$$
 corresponds to $s_\ast^1:
(X,T_1)\to (Y,T_2)$.

Check now compatibility of the transition $s_\ast^1\to s_\ast^2$
with the multiplication of morphisms. Given $s_1: W(X)\to W(Y)$
and $s_2: W(Y)\to W(Z)$, we have $s_2s_1: W(X)\to W(Z)$. Using
once more the fact that the transition $\Theta^0\to
\Hal_{\Theta}(\Phi)$ is a functor, we get
$$
(s_2^1s_1^1)_\ast= s_{2\ast}^1s_{1\ast}^1,
$$
$$
(s_2^2s_1^2)_\ast= s_{2\ast}^2s_{1\ast}^2,
$$

Apply $\widetilde\beta$. We need to verify that
$\widetilde\beta(s_{2\ast}^1s_{1\ast}^1)=
\widetilde\beta(s_{2\ast}^1)\widetilde\beta(s_{1\ast}^1). $ We
have
$$
\widetilde\beta(s_{2\ast}^1s_{1\ast}^1)=\widetilde\beta(s_{2}^1s_{1}^1)_\ast=
(s_{2}^2s_{1}^2)_\ast=
s_{2\ast}^2s_{1\ast}^2=\widetilde\beta(s_{2\ast}^1)\widetilde\beta(s_{1\ast}^1).
$$
This gives compatibility with the multiplication as well as with
the unit. Hence, we have the functor $\widetilde\beta:
L_\Theta(\Phi_1)\to L_\Theta(\Phi_2)$.

\subhead{ 5.4 Knowledge bases}
\endsubhead

We proceed from a multi-model $(G, \Phi, F)$. A multi-model  $(G,
\Phi, F)$ defines a system of models  $(G, \Phi, f,)$ where $f$
runs the set $F$. Here $G$ is an algebra in $\Theta$, and $\Phi$
is a set of relations. Recall that both the algebra $G\in \Theta$
and a relation $f\in F$ are multi-sorted.  The set $F$ is a set of
instances $f$, where $f$ is a interpretation of the set $\Phi$ in
$G$.

To every such multi-model corresponds  a knowledge base $KB =
KB(G,\Phi,F)$. The definition slightly differs from that of
\cite{PTP}.

\proclaim{Definition} A knowledge base $KB = KB(G,\Phi,F)$
consists of two categories. The first one is the category of
knowledge description $L_\Theta(\Phi)$, and the second one is the
category of knowledge content $K_{\Phi\Theta}(f)$. These two
categories are related by the functor
$$
\Ct_f : L_{\Theta}(\Phi) \to K_{\Phi\Theta}(f).
$$
\endproclaim

This functor $\Ct_f$ transforms knowledge description to content
of knowledge.
We do not assume that between different $f_1$ and $f_2$ in $F$
there are any ties: instances are independent. On the other hand,
between some $f_1$ and $f_2$ there may be relations that we will
try to take into account (see Section 7).

A content of knowledge $\Ct_f(X,T)=(X,T^f)$ corresponds to an
object $(X,T)$ of the category $L_\Theta(\Phi)$, which is a
description of knowledge. We view the description $T$ as a {\it
query} to a knowledge base, and $A=T^f$ as a {\it reply to this
query}.

Besides, if there is a relation $s_\ast$ between $(X,T_1)$ and
$(Y,T_2)$, then  there will be a relation $\widetilde s
=\widetilde {s_\ast}$ between $(X,A)$ and $(Y,B)$, where $A=T_1^f,
\ B=T_2^f$  .

This peculiarity of the definition naturally reflects geometrical
essence of knowledge.

In fact, in this definition of a knowledge base the category of
knowledge is decomposed to two categories: the category of
description of knowledge and the category of content of knowledge,
tied by the functor of transition from description to content.

\head{ 6. Equivalence of  knowledge bases}
\endhead

\subhead{ 6.1 Definition}
\endsubhead

Let the knowledge bases $KB_1=KB(G_1,\Phi_1,F_1)$  and
$KB_2=KB(G_2,\Phi_2,F_2)$ correspond to the given multi-models
$(G_1,\Phi_1,F_1)$ and $(G_2,\Phi_2,F_2)$.

\proclaim{Definition 1} Knowledge bases $KB_1$ and $KB_2$ are
called informationally  equivalent, if there exists a bijection
$\alpha: F_1\to F_2$ such that for every $f\in F_1$ there exist
homomorphisms
$$
\beta_f: \Hal_\Theta(\Phi_1)\to \Hal_\Theta(\Phi_2)
$$
$$
\beta_f': \Hal_\Theta(\Phi_2)\to \Hal_\Theta(\Phi_1)
$$
and an isomorphism of categories
$$
\widetilde\gamma_f: K_{\Phi_1\Theta}(f)\to
K_{\Phi_2\Theta}(f^\alpha)
$$
such that the commutative diagrams of functors of categories hold:
\endproclaim
$$
\CD
 L_\Theta(\Phi_1)@>\widetilde\beta_f>> L_\Theta(\Phi_2)\\
@V\Ct_f VV @VV \Ct_{f^\alpha}V\\
K_{\Phi_1\Theta}(f)@>\widetilde\gamma_f>> K_{\Phi_2\Theta}(f^\alpha)\\
\endCD
$$
and
$$
\CD
 L_\Theta(\Phi_1)@<\widetilde\beta'_f<< L_\Theta(\Phi_2)\\
@V\Ct_f VV @VV \Ct_{f^\alpha}V\\
K_{\Phi_1\Theta}(f)@<(\widetilde\gamma_f)^{-1}<< K_{\Phi_2\Theta}(f^\alpha)\\
\endCD
$$

Denote these diagrams by $\ast$ and $\ast\ast$ respectively.
Rewrite commutative diagrams for the object $(X, T)$ of the
category $L_\Theta (\Phi_1)$ in the form $(X, T^f)^{\widetilde
\gamma _f} = (X, T^{\beta_f f^\alpha})$ and for the object $(X,
T)$ of the category $L_\Theta (\Phi_2)$ in the form $(X,
T^{f^\alpha})^{\widetilde \gamma _f ^{-1}} = (X, T^{\beta' _f
f})$.

From this follows
$$
(X, T^f) = (X, T^{\beta_f
f^\alpha})^{\widetilde{(\gamma_f)}^{-1}},
$$
$$
(X, T^{f^\alpha}) = (X, T^{\beta'_f f})^{\widetilde \gamma_f }.
$$

The last means that everything which can be known from $KB_1$ can
be also known from $KB_2$ and vice versa. Similar property holds
for morphisms, i.e. for relations between objects. Equivalence of
knowledge bases we consider as a triple $(\alpha,\ast,\ast\ast)$,
where $\alpha: F_1\to F_2$ is a bijection, while $\ast$ and
$\ast$$\ast$ define the corresponding diagrams for every $f\in
F_1$.

The next proposition deals with the transition from  knowledge
bases to data\-ba\-ses. Let $R_f$ be the image of the homomorphism
$\Val_f: \Hal_\Theta(\Phi)\to \Hal_\Theta(G).$

\proclaim{Proposition 3} If a bijection $\alpha: F_1\to F_2$
determines equivalence of the bases $KB_1$ and $KB_2$ then for
every $f\in F_1$ we have an isomorphism of Halmos algebras
$\gamma_f:R_f\to R_{f^\alpha}$.
\endproclaim

Proof.

Proceed from the corresponding diagrams $\ast$ and $\ast\ast$.
Given  a set $X$, take a set $T$ consisting of one element
$u\in\Hal_{\Phi_1\Theta}(X)$. In this case $T^f=\Val_f(u)$. We
have $Ct_f(X,T)=(X,\Val_f(u))$,
$$
(X,\Val_f(u))^{\widetilde\gamma_f}=Ct_{f^\alpha}(X,u^\beta)=(X,(u^\beta)^{f^\alpha})=
(X,\Val_{f^\alpha}(u^\beta)).
$$

Hence, $\widetilde \gamma_f$ transfers $\Val_f(u)$ to
$\Val_{f^\alpha}(u^\beta)$ for every $u$, which means that
$\widetilde\gamma_f$ induces a mapping $\gamma_f: R_f\to
R_{f^\alpha}$. It is a homomorphism since $\Val_f$ and $\beta$ are
homomorphisms of algebras, and it is an injection since every
$R_f$ is a simple algebra \cite{Pl1}.

Let now $u_1$ be an arbitrary element of $\Hal_{\Phi_2\Theta}(X).$
Then the second diagram gives
$$
(X,\Val_{f^\alpha}(u_1))^{\widetilde\gamma_f^{-1}}=
(X,\Val_f(u_1^{\beta'_f})),
$$
and
$$
(X,\Val_{f^\alpha}(u_1))=(X,\Val_f(u_1^{\beta'_f}))^{\widetilde\gamma_f}=
(X,\Val_{f^\alpha}(u))^{\widetilde\gamma_f},
$$
where $u=u_1^{\beta'_f}$. This implies that $\gamma_f:R_f\to
R_{f^\alpha}$ is a surjection. Hence, we have an isomorphism
$\gamma_f:R_f\to R_{f^\alpha}$.


\subhead{7.2 Finite models}
\endsubhead

First of all it is clear that for finite models $(G,\Phi, F)$ the
corresponding KB remains, in general, infinite.

We prove the following main

\proclaim{Theorem 4} Let the given models be finite. Then the
knowledge bases $KB_1$ and $KB_2$ are equivalent if and only if
there exists a bijection $\alpha: F_1\to F_2$ such that for every
$f\in F_1$ there is an isomorphism $\gamma_f: R_f\to
R_{f^\alpha}$.
\endproclaim

Proof.

In one direction the statement is always true. Let  now $\gamma_f:
R_f\to R_{f^\alpha}$ be an isomorphism for every $f\in F_1$.
According to  Theorem 4 from   \cite{PT}  there are
 the
homomorphisms $\beta_f:\Hal_{\Theta}(\Phi_1)\to
\Hal_{\Theta}(\Phi_2)$ and $\beta'_f:\Hal_{\Theta}(\Phi_2)\to
\Hal_{\Theta}(\Phi_1)$ such that the diagrams

$$
\CD
 \Hal_{\Theta}(\Phi_1)@>\beta_f>> \Hal_{\Theta}(\Phi_2)\\
@V\Val_f VV @VV\Val_{f^\alpha} V\\
R_f@>\gamma_f>> R_{f^\alpha}\\
\endCD
$$

$$
\CD
 \Hal_{\Theta}(\Phi_1)@<\beta'_f<< \Hal_{\Theta}(\Phi_2)\\
@V\Val_f VV @VV\Val_{f^\alpha} V\\
R_f@<\gamma_f^{-1}<< R_{f^\alpha}\\
\endCD
$$

\noindent
 are commutative.

 Simultaneously, there are functors
 $$
\widetilde\beta_f: L_{\Theta}(\Phi_1)\to L_{\Theta}(\Phi_2),
$$
$$
\widetilde\beta'_f: L_{\Theta}(\Phi_2)\to L_{\Theta}(\Phi_1).
$$

 It is left to define the isomorphism of categories $
\widetilde{\gamma_f}:K_{\Phi_1\Theta}(f)\to
K_{\Phi_2\Theta}(f^\alpha)$ such that the diagrams of the types
$\ast$ and $\ast\ast$ be commutative.

First we define $\widetilde\gamma_f$ on objects and then on
morphisms. Take an object (X,T) of the category $L_\Theta(\Phi_1)$
for an arbitrary object $(X,A)$ of the category
$K_{\Phi_1\Theta}(f)$ with $T^f=A$. We have
$Ct_f(X,T)=(X,T^f)=(X,A)$. Set
$$
(X,A)^{\widetilde\gamma_f}=(X,T^f)^{\widetilde\gamma_f}= (X,
\bigcap_{u\in T}\gamma_f\Val_f(u))= $$
$$(X, \bigcap_{u\in
T}\Val_{f^\alpha}(u^{\beta_f}))=(X,T^{\beta_ff^\alpha}).
$$

We want to show that this definition does not depend on the choice
of the set $T$ with $T^f=A$. Consider first the case when
$T_1^f=T_2^f=A$ and the sets $T_1$ and $T_2$ are finite. We have:
$(X,T_1^f)^{\widetilde\gamma_f}= (X,T_1^{\beta_ff^\alpha})$ and
$(X,T_2^f)^{\widetilde\gamma_f}= (X,T_2^{\beta_ff^\alpha}).$

We need to check that $T_1^{\beta_ff^\alpha}=
T_2^{\beta_ff^\alpha}. $ Indeed,
$$
T_1^{\beta_ff^\alpha}=\bigcap_{u_1\in
T_1}\Val_{f^\alpha}(\beta_fu_1))=\bigcap_{u_1\in
T_1}\gamma_f\Val_f(u_1)).
$$
Since $\gamma_f:R_f\to R_{f^\alpha}$ is an isomorphism of algebras
and $T_1,$ $T_2$ are finite sets, we can rewrite the expression in
the form
$$
T_1^{\beta_ff^\alpha}=\gamma_f(\bigcap_{u_1\in
T_1}\Val_{f}(u_1))=\gamma_f(\bigcap_{u_2\in T_2}\Val_f(u_2))=
\bigcap_{u_2\in T_2}\gamma_f\Val_{f}(u_2))=T_2^{\beta_ff^\alpha}.
$$

Passing to the general case we proceed from finite models. Every
finite model is geometrically noetherian, i.e., if
$A=T_1^f=T_2^f$, then in $T_1$ and $T_2$ one can find finite
subsets $T_{01}$ and $T_{02}$ with $T{_{01}}^f=T_{02}^f=A$. Here,
$T_{01}^{\beta_ff_\alpha}=T_{02}^{\beta_ff_\alpha}$. We have to
verify that $ T_1^{\beta_ff^\alpha}=T_2^{\beta_ff^\alpha} $ and
$T_{01}^{\beta_ff_\alpha}=\bigcap_{u_1\in
T_1}\Val_{f^\alpha}(\beta_fu_1)$. We can take a finite subset
$T_{10}$ in $T_1$ such that $ T_1^{\beta_ff^\alpha}=
T_{10}^{\beta_ff^\alpha}$. Take the union of sets $T_{10}$ and
$T_{01}$ and denote it by $T_{001}$. Then $T_{001}^f=A=T_1^f$,
$T_1^{\beta_ff^\alpha}=T_{001}^{\beta_ff^\alpha}$. Analogously,
for $T_2$ take $T_{002}$ and $A=T_{001}^f=T_{002}^f$. Besides
that,
$$
T_1^{\beta_ff^\alpha}=T_{001}^{\beta_ff^\alpha}=
T_2^{\beta_ff^\alpha}.
$$
The equality $ T_1^{\beta_ff^\alpha}= T_2^{\beta_ff^\alpha}$ gives
commutativity of the diagram for objects.

Similarly, we build $\widetilde\gamma_f^{-1}$ having
$\gamma_f^{-1}$ and the equality
$\widetilde\gamma_f^{-1}=\widetilde{\gamma_f^{-1}}$ holds.

Now let us pass to morphisms. Remind first of all that to every
homomorphism $s:W(Y)\to W(X) $ there correspond
$$
s_\ast^1:\Hal_{\Phi_1\Theta}(Y) \to \Hal_{\Phi_1\Theta}(X),
$$
$$
s_\ast^2: \Hal_{\Phi_2\Theta}(Y) \to \Hal_{\Phi_2\Theta}(X).
$$

Let the objects $(Y,T_2)$ and $(X,T_1)$ be given in
$L_\Theta(\Phi_1)$. Recall that $s$ is admissible for $T_2$ and
$T_1$ if $ s_\ast^1(u)\in T_1$ for every $u\in T_2$. Here
$s_\ast^1:(Y,T_2)\to (X,T_1)$ is a morphism. Proceed further from
an arbitrary homomorphism
$\beta:\Hal_\Theta(\Phi_1)\to\Hal_\Theta(\Phi_2)$. It had been
proved that if $s$ is admissible for $T_2$ and $T_1$ then the same
$s$ is admissible for $T_2^\beta$ and $T_1^\beta$ as well, i.e., $
s_\ast^1(u)\in T_1^\beta$ for every $u\in T_2^\beta$. Hence, we
have a morphism
$$
\widetilde\beta(s_\ast^1)=s_\ast^2:(Y,T_2^\beta)\to(X,T_1^\beta).
$$
Take now $\beta=\beta_f$ and apply $\Ct_{f^\alpha}$:
$$
\Ct_{f^\alpha}(s^2_\ast):\Ct_{f^\alpha}(X,T_1^{\beta_X})\to
\Ct_{f^\alpha}(Y, T_2^{\beta_X}).
$$
It can be rewritten as
$$
\Ct_{f^\alpha}(s_\ast^2):(X,T_1^{\beta_Xf^\alpha})\to (X,
T_2^{\beta_Yf^\alpha})
$$
or $$ \Ct_{f^\alpha}(s_\ast^2):(X,T_1^f)^{\gamma_f}\to (Y,
T_2^{f})^{\gamma_f}.
$$
Let now $T_1^f=A$, $T_2^f=B$. For $s_\ast^1:(Y,T_2)\to (X,T_1)$ we
have
$$
\Ct_{f}(s_\ast^1):(X,T_1^f)\to (Y, T_2^f)
$$
and  a related morphism
$$
\Ct_{f^\alpha}(s_\ast^2):(X,T_1^f)^{\gamma_f}\to (Y,
T_2^{f_2})^{\gamma_f}.
$$

Commutativity of the diagram on morphisms means that
$$
\widetilde\gamma_f\Ct_f(s_\ast^1)=\Ct_{f^\alpha}(\widetilde\beta_f(s_\ast^1))
$$
for every $s^1_\ast: (Y,T_2)\to (X,T_1). $

Continuing consideration of finite models, proceed from the
isomorphism $\gamma_f:R_f\to R_{f^\alpha}$ and the corresponding
functor $\widetilde\gamma_f: K_{\Phi_1\Theta}(f)\to
K_{\Phi_2\Theta}(f^\alpha)$. This functor had been defined on the
objects, and now we are going to define it on morphisms.

Let $\tau: (X,A)\to (Y,B)$ be a morphism in $K_{\Phi_1\Theta}(f)$.
This $\tau$ appears as follows. A morphism
$$
s_\ast^1:\Hal_{\Phi_1\Theta}(Y) \to \Hal_{\Phi_1\Theta}(X)
$$
corresponds to $s: W(Y)\to W(X)$. If now $A=T_1^f$, $B=T_2^f$ and
$s^1_\ast$ is admissible for  $T_2$ and $T_1$ then we have
$\widetilde s_\ast^1: (X,A)\to (Y,B). $ We may say that
 $\tau=\widetilde s_\ast^1$ for some $s_\ast^1$.

 Define
 $$
 \widetilde\gamma_f(\widetilde s_\ast^1)
 =\widetilde s_\ast^2: (X,
 T_1^f)^{\widetilde\gamma_f}\to (Y,T_2^f)^{\widetilde\gamma_f}.
 $$
Here,
$$
(X, T_1^f)^{\widetilde\gamma_f}=(X,T_1^{\beta_ff^\alpha}),
 $$
  $$
 (Y, T_2^f)^{\widetilde\gamma_f}=(Y,T_2^{\beta_ff^\alpha})
$$
do not depend on the choice of $T_1$ and $T_2$ with $T_1^f=A$ and
$T_2^f=B$. Check further that $\widetilde\gamma_f:
K_{\Phi_1\Theta}(f)\to K_{\Phi_2\Theta}(f^\alpha)$ determined in
such a way is in fact a functor and this functor provides
commutativity of the diagram on morphisms.

Note first of all that the definition of $\widetilde\gamma_f$ on
morphisms can be rewritten as
$$
\widetilde\gamma_f(\Ct_f(s_\ast^1))=(\Ct_{f^\alpha}(s_\ast^2)).
$$
Take two morphisms $\widetilde{s_{1\ast}^1}=\Ct_f(s_{1\ast}^1)$
and $\widetilde{s_{2\ast}^1}=\Ct_f(s_{2\ast}^1)$ and consider the
product
$$
\widetilde{s_{1\ast}^1}\widetilde{s_{2\ast}^1}=\Ct_f(s_{1\ast}^1)\Ct_f(s_{2\ast}^1)=
\Ct_f(s_{2\ast}^1s_{1\ast}^1)=\widetilde{s_{2\ast}^1s_{1\ast}^1}=
\widetilde{(s_2s_1)_\ast^1}.
$$

Apply $\widetilde\gamma_f$:
$$
\widetilde\gamma_f((\widetilde{s_2s_1})_\ast^1)=(\widetilde{(s_2s_1)_\ast^2})=
\widetilde{s_{2\ast}^2s_{1\ast}^2}=\widetilde{s_{1\ast}^2}\widetilde{s_{2\ast}^2}=
\widetilde\gamma_f(\widetilde{s_{1\ast}^1})\widetilde\gamma_f(\widetilde{s_{2\ast}^1}).
$$
Now check the commutativity of the diagram
$$
\CD
 L_{\Theta}(\Phi_1)@>\widetilde\beta_X>> L_{\Theta}(\Phi_2)\\
@V\Ct_f VV @VV\Ct_{f^\alpha} V\\
K_{\Phi_1\Theta}(f)@>\widetilde\gamma_f>> K_{\Phi_1\Theta}(f^\alpha)\\
\endCD
$$

Take a morphism $s_\ast^1:(Y,T_2)\to(X,T_1)$ in
$L_\Theta(\Phi_1)$. We have
$$
\widetilde{\beta_X}
(s_\ast^1):(Y,T_2^{\beta_X})\to(X,T_1^{\beta_X}),
$$
and
$$
\Ct_{f^\alpha}\widetilde{\beta_X}
(s_\ast^1):(X,T_1^{\beta_Xf^\alpha})\to(Y,T_2^{\beta_Xf^\alpha}).
$$
Rewrite it as
$$
\Ct_{f^\alpha}\widetilde{\beta_X}
(s_\ast^1):(X,T_1^f)^{\widetilde\gamma_f}\to(Y,T_2^{f})^{\widetilde\gamma_f}.
$$
Further,
$$
\Ct_{f} (s_\ast^1):(X,T_1^f)\to(Y,T_2^{f}),
$$
$$
\widetilde\gamma_f\Ct_{f}
(s_\ast^1):(X,T_1^f)^{\widetilde\gamma_f}\to(Y,T_2^{f})^{\widetilde\gamma_f}.
$$
Check now the equality
$$
\widetilde\gamma_f\Ct_{f}
(s_\ast^1)=\Ct_{f^\alpha}\widetilde\beta_X (s_\ast^1)
$$
for every $ s_\ast^1$. We have
$$
\widetilde\gamma_f\Ct_{f} (s_\ast^1)=\widetilde\gamma_f(\widetilde
s_\ast^1)=\widetilde s_\ast^2,
$$
$$
\Ct_{f^\alpha}\widetilde\beta_X(s_\ast^1)=\Ct_{f^\alpha}(s_\ast^2)=
\widetilde s_\ast^2.
$$
This gives commutativity of the diagram $\ast$ of morphisms, i.e.,
$$
\widetilde\gamma_f\Ct_{f} =\Ct_{f^\alpha}\widetilde\beta_X.
$$
The same can be done for the functor
$\widetilde{\gamma_f^{-1}}=\widetilde\gamma_f^{-1}$ and the second
commutative diagram $\ast\ast$ that finishes the proof of the
theorem

\subhead{7. Additional remarks}
\endsubhead

7.1. Let us look at the definition of equivalence from the general
perspective of category theory. Given two functors $\varphi_1 :
C_1 \to C_1 ^0$ and $\varphi_2 : C_2 \to C_2 ^0$, we say that
$C_1$ and $C_2$ are equivalent in respect to $\varphi_1$ and
$\varphi_2$, if there is an isomorphism $\psi : C_1 ^0 \to C_2 ^0$
and functors $\psi_1 : C_1 \to C_2$, $\psi_2 : C_2 \to C_1$ with
the commutative diagrams
$$
\CD
 C_1@>\psi_1>> C_2\\
@V\varphi_1 VV @VV\varphi_2 V\\
C_1^0@>\psi>> C_2^0\\
\endCD
$$

$$
\CD
 C_1@<\psi_2<< C_2\\
@V\varphi_1 VV @VV\varphi_2 V\\
C_1^0@<\psi^{-1}<< C_2^0\\
\endCD
$$

Usual equivalence of categories is equivalence in respect to the
transition to skeletons of categories. In our situation we may say
that equivalence of knowledge bases means that there exists
equivalence of categories of description of knowledge in respect
to transition to the categories of knowledge content.

7.2. Let us return to the definition of knowledge bases with
multi-models $(G_1, \Phi_1,$ $ F_1)$ and $(G_2, \Phi_2, F_2)$, and
let the bijection $\alpha : F_1 \to F_2$ determine equivalence of
the corresponding $KB_1$ and $KB_2$. Assume that two instances
$f_1$ and $f_2$ from $F_1$ are connected by a commutative diagram
$$
\CD
\Hal_{\Theta}(\Phi_1) @>\Val_{f_1}>> R_{f_1} \\
@. @/SE/ \Val_{f_2}// @VV\gamma V\\
@. R_{f_2}\\
\endCD
$$
where $\gamma$ is a homomorphism of algebras. We want to evaluate
the relation between $f_1 ^\alpha$ and $f_2 ^\alpha$.

Proceed from the diagrams

$$
\CD
 \Hal_{\Phi_1\Theta}@>\beta_f>> \Hal_{\Phi_2\Theta}\\
@V\Val_f VV @VV\Val_{f^\alpha} V\\
R_f@>\gamma_f>> R_{f^\alpha}\\
\endCD
$$

$$
\CD
 \Hal_{\Phi_1\Theta}@<\beta' _f<< \Hal_{\Phi_2\Theta}\\
@V\Val_f VV @VV\Val_{f^\alpha} V\\
R_f@<\gamma_f^{-1}<< R_{f^\alpha}\\
\endCD
$$

$$
\CD
 R_{f_1}@>\gamma_{f_1}>> R_{f_1^\alpha}\\
@V\gamma VV @VV\gamma^\alpha V\\
R_{f_2}@>\gamma_{f_2}>> R_{f_2^\alpha}\\
\endCD
$$

Here,
$$
\gamma\Val_{f^1}=\Val_{f^2}, \qquad
\gamma^\alpha=\gamma_{f_2}\gamma\gamma_{f_1}^{-1}
$$
and
$$
\gamma^\alpha\Val_{f_1^\alpha}=\gamma^\alpha\gamma_{f_1}\Val_{f_1}\beta'_{f_1}=
\gamma_{f_2}\gamma\Val_{f_1}\beta'_{f_1}=\gamma_{f_2}\Val_{f_2}\beta'_{f_1}=
\Val_{f_2^\alpha}\beta_{f_2}\beta'_{f_1}.
$$
Hence,
$\gamma^\alpha\Val_{{f^1}^\alpha}=\Val_{f_2^\alpha}\beta_{f_2}\beta'_{f_1},$
i.e., the connection is twisted by the product
$\beta_{f_2}\beta'_{f_1}.$

 At last, let us note that from the
diagrams above follow the natural identities:

1. $\Val_{f}(u)=\Val_{f}(\beta'_{f}\beta_{f}(u))$ for every
$u\in\Hal_\Theta(\Phi_1)$.

2. $\Val_{f^\alpha}(u)=\Val_{f^\alpha}(\beta_{f}\beta'_{f}(u)$ for
every $u\in\Hal_\Theta(\Phi_2)$.

7.3. Note that the equivalence condition of two knowledge bases in
the case of finite multi-models can be formulated in terms of
these multi-models (cf. \cite{PTP}).

\proclaim{Definition 2} Let the models $(G_1,\Phi_1,f_1)$ and
$(G_2,\Phi_2,f_2)$ be given. Let $\Aut(f_1)$ and  $\Aut(f_2)$ be
the corresponding groups of automorphisms. The models
$(G_1,\Phi_1,f_1)$  and $(G_2,\Phi_2,f_2)$ are called automorphic
equivalent if there exists an isomorphism of algebras $\delta:
G_1\to G_2$ such that
 $$
 \Aut(f_2)=\delta \Aut(f_1)\delta^{-1}.
 $$
\endproclaim
\proclaim{Definition 3} Let the multi-models $(G_1,\Phi_1,F_1)$
and $(G_2,\Phi_2,F_2)$ be given. These multi-models are called
automorphic equivalent if there exists a bijection $\alpha :
F_1\to F_2$ such that for every $f\in F_1$ the models
$(G_1,\Phi_1,f)$  and $(G_2,\Phi_2,f^\alpha)$ are automorphic
equivalent.
\endproclaim

It is natural to define an isomorphism of multi-models with the
same set of relations $\Phi_1$ and $\Phi_2$. An isomorphism of
multi-models implies their automorphic equivalence. Evidently, the
inverse statement is not true.

Let the knowledge bases $KB_1=KB(G_1,\Phi_1,F_1)$ and
$KB_2=KB(G_2,\Phi_2,F_2)$ with the finite multi-models be given.

\proclaim{Theorem 5}The knowledge bases $KB_1=KB(G_1,\Phi_1,F_1)$
and $KB_2=KB(G_2,$ $\Phi_2,F_2)$ are informationally equivalent if
and only if the corresponding models are automorphic equivalent.
\endproclaim

The proof of this theorem is parallel to the proof of the
corresponding theorem in \cite{PTP} and uses the Galois-Krasner
theory in the given variety of algebras $\Theta$ \cite{Pl1}.
Theorem 5 provides an algorithm for the informational equivalence
verification.

 \bye

\ref \key {\bf [JCP]} \by Jeavons, P., Cohen, D.,; Pearson, J.
\paper Constraints and universal algebra \jour Ann. Math.
Artificial Intelligence \vol 24\year 1999\pages 51--67

\heading{Conclusion}
\endheading

We can speak about knowledge and a system of knowledge. As  a
rule, considerations of knowledge deal with a fixed domain of
knowledge or of a system of knowledge. These domains could be, for
example, political, historical or archeological knowledge,
humanity and scientific knowledge, knowledge in the field of
mathematics, or, more specifically, geometry, logic, etc.
Knowledge in one field can be based on the knowledge in other
fields, on experiments, graphs, visual information.

In the science of knowledge we consider only knowledge which admit
formalization in some logic. The logic may be different. It is
often oriented towards the corresponding field of knowledge.


In the paper we consider a special situation. We confine ourselves
with the First Order Logic (FOL), but along with logic we use also
geometry. We define formally the notion of knowledge and study the
category of knowledge. We give definition of a mathematical model
of a knowledge base. This definition illuminates a clear
distinction between knowledge base and database concepts.

The problem of informational equivalence of knowledge bases has
not been stated earlier, at least in  a way that allows to speak
of algorithm of verification of informational equivalence.
Formerly only the problem of informational equivalence of
databases had been treated \cite{12}.

 Informational equivalence of knowledge bases is a
very difficult problem that till now had no formal solution.
Mathematical model allows to define exactly the notion of
informational equivalence of knowledge bases and construct an
algorithm of verification of knowledge bases equivalence.

\bye